\documentclass[11pt]{article}
\usepackage{amssymb}
\usepackage[dvips]{epsfig}
\usepackage{times}
\usepackage{graphicx}

\newtheorem{Lemma}{Lemma}[section]
\newtheorem{Theorem}{Theorem}
\newtheorem{Proposition}[Lemma]{Proposition}
\newtheorem{Corollary}[Lemma]{Corollary}
\newtheorem{Remark}[Lemma]{Remark}
\newtheorem{Definition}[Lemma]{Definition}

\newenvironment{Proof}%
 {\begin{trivlist} \item[]{\bf Proof. }}%
 {\hspace*{\fill}$\rule{.3\baselineskip}{.35\baselineskip}$\end{trivlist}}
 \newenvironment{Proof1}%
{\begin{trivlist} \item[]{\bf Proof }}%
{\hspace*{\fill}$\rule{.3\baselineskip}{.35\baselineskip}$\end{trivlist}}
 {\begin{trivlist}\item[]\textbf{Acknowledgments }}{\end{trivlist}}

\makeatletter \@addtoreset{equation}{section} \makeatother

\newcommand{\C}{\mathbb{C}}

\newcommand{\R}{\mathbb{R}}

\def\Im{\mathop{\mathrm{Im}}}

\newfam\bifam
\font\tenbi=cmmib10 scaled \magstep1 \font\sevenbi=cmmib10 at 11pt
\font\fivebi=cmmib10 at 6pt \textfont\bifam = \tenbi
\scriptfont\bifam= \sevenbi \scriptscriptfont\bifam= \fivebi

\begin{document}

\title{\bf Count of eigenvalues in the generalized eigenvalue problem}

\author{Marina Chugunova and Dmitry Pelinovsky  \\
{\small $^{1}$ Department of Mathematics, McMaster University,
Hamilton, Ontario, Canada, L8S 4K1} }

\date{\today}
\maketitle

\begin{abstract}
We address the count of isolated and embedded eigenvalues in a
generalized eigenvalue problem defined by two self-adjoint operators
with a positive essential spectrum and a finite number of isolated
eigenvalues. The generalized eigenvalue problem determines spectral
stability of nonlinear waves in a Hamiltonian dynamical system. The
theory is based on the Pontryagin's Invariant Subspace theorem in an
indefinite inner product space but it extends beyond the scope of
earlier papers of Pontryagin, Krein, Grillakis, and others. Our main
results are (i) the number of unstable and potentially unstable
eigenvalues {\em equals} the number of negative eigenvalues of the
self-adjoint operators, (ii) the total number of isolated
eigenvalues of the generalized eigenvalue problem is {\em bounded
from above} by the total number of isolated eigenvalues of the
self-adjoint operators, and (iii) the quadratic form defined by the
indefinite inner product is strictly positive on the subspace
related to the absolutely continuous part of the spectrum of the
generalized eigenvalue problem. Applications to solitons and
vortices of the nonlinear Schr\"{o}dinger equations and solitons of
the Korteweg--De Vries equations are developed from the general
theory.
\end{abstract}


\newpage

\section{Introduction}

Stability of equilibrium points in a Hamiltonian system of finitely
many interacting particles is defined by the eigenvalues of the
generalized eigenvalue problem \cite{GKrein},
\begin{equation}
\label{eigenvalue-finite-dimension} A {\bf u} = \gamma K {\bf u},
\qquad {\bf u} \in \mathbb{R}^n,
\end{equation}
where $A$ and $K$ are symmetric matrices in $\mathbb{R}^{n \times 
n}$ which define the quadratic forms for potential and kinetic 
energies, respectively. The eigenvalue $\gamma$ corresponds to the 
normal frequency $\lambda = i \omega$ of the normal mode of the 
linearized Hamiltonian system near the equilibrium point, such that 
$\gamma = - \lambda^2 = \omega^2$. The linearized Hamiltonian system 
is said to have an unstable eigenvalue $\gamma$ if $\gamma < 0$ or 
$\Im(\gamma) \neq 0$. Otherwise, the system is weakly spectrally 
stable. Moreover, the equilibrium point is a minimizer of the 
Hamiltonian if all eigenvalues $\gamma$ are positive and semi-simple 
and the quadratic forms for potential and kinetic energies evaluated 
at eigenvectors of (\ref{eigenvalue-finite-dimension}) are strictly 
positive.

When the matrix $K$ is positive definite, all eigenvalues $\gamma$
are real and semi-simple (that is the geometric and algebraic
multiplicities coincide). By the Sylvester's Inertia Law theorem
\cite{Gelfand}, the numbers of positive, zero and negative
eigenvalues of the generalized eigenvalue problem
(\ref{eigenvalue-finite-dimension}) {\em equal to} the numbers of
positive, zero and negative eigenvalues of the matrix $A$. When $K$
is not positive definite, a complete classification of eigenvalues
$\gamma$ in terms of real eigenvalues of $A$ and $K$ has been
developed with the use of the Pontryagin's Invariant Subspace
theorem \cite{Pontryagin}, which generalizes the Sylvester' Inertia
Law theorem.

We are concerned with spectral stability of spatially localized
solutions in a Hamiltonian infinite-dimensional dynamical system. In
many problems, a linearization of the nonlinear system at the
spatially localized solution results in the generalized eigenvalue
problem of the form (\ref{eigenvalue-finite-dimension}) but $A$ and
$K^{-1}$ are now self-adjoint operators on a complete
infinite-dimensional metric space. There has been recently a rapidly
growing sequence of publications on mathematical analysis of the
spectral stability problem in the context of nonlinear
Schr\"{o}dinger equations and other nonlinear evolution equations
\cite{CPV,GKP,KKS,KP,KSchlag,Pel}. Besides predictions of spectral
stability or instability of spatially localized solutions in
Hamiltonian dynamical systems, linearized Hamiltonian systems are  
important in analysis of orbital stability \cite{GSS1,GSS2,CP03}, 
asymptotic stability \cite{Per04,RSS,Cuc}, stable manifolds 
\cite{CCP05,Schlag}, and blow-up of solutions in nonlinear equations 
\cite{Per}.

It is the purpose of this article to develop analysis of the
generalized eigenvalue problem in infinite dimensions by using the
Pontryagin space decomposition \cite{Pontryagin}. The theory of
Pontryagin spaces was developed by M.D. Krein and his students (see
books \cite{Azizov,krein,iohvidov}) and partly used in the context
of spectral stability of solitary waves by MacKay \cite{MacKay},
Grillakis \cite{G90}, and Buslaev \& Perelman \cite{BusPer} (see
also a recent application in \cite{GKP}). We shall give an elegant
geometric proof of the Pontryagin's Invariant Subspace theorem. An
application of the theorem recovers the main results obtained in
\cite{CPV,KKS,Pel}. Moreover, we obtain a new {\em inequality} on
the number of positive eigenvalues of the linearized Hamiltonian
that extends the count of eigenvalues of the generalized eigenvalue 
problem.

The structure of the paper is as follows. Main formalism of the
generalized eigenvalue problem is described in Section 2. The
Pontryagin's Invariant Subspace theorem is proved in Section 3. Main
results on eigenvalues of the generalized eigenvalue problem are
formulated and proved in Section 4. Section 5 contains applications
of the main results to solitons and vortices of the nonlinear
Schrodinger equations and solitons of the Korteweg--De Vries
equations.

\section{Formalism}

Let $L_+$ and $L_-$ be two self-adjoint operators defined on the
Hilbert space ${\cal X}$ with the inner product $(\cdot,\cdot)$. Our
main assumptions are listed below.

\begin{itemize}
\item[P1] The essential spectrum of $L_{\pm}$ in ${\cal X}$
includes the absolute continuous part $\sigma_c(L_{\pm})$ bounded
from below by $\omega_+ \geq 0$ and $\omega_- > 0$ and finitely many
embedded eigenvalues of finite multiplicities.

\item[P2] The discrete spectrum of $L_{\pm}$ in ${\cal X}$
includes finitely many isolated eigenvalues of finite multiplicities
with $p(L_{\pm})$ positive, $z(L_{\pm})$ zero, and $n(L_{\pm})$
negative eigenvalues.
\end{itemize}

We shall consider the linearized Hamiltonian problem defined by the
self-adjoint operators $L_{\pm}$ in ${\cal X}$,
\begin{equation}
\label{coupled-problem} L_+ u = - \lambda w, \qquad L_- w = \lambda
u,
\end{equation}
where $\lambda \in \mathbb{C}$ and $(u,w) \in {\cal X} \times {\cal
X}$. By assumption P1, the kernel of $L_-$ is isolated from the
essential spectrum. Let ${\cal H}$ be the constrained Hilbert space,
\begin{equation}
\label{constraint} \mathcal{H} = \left\{ u \in {\cal X} : \;\; u
\perp {\rm ker}(L_-) \right\},
\end{equation}
and let ${\cal P}$ be the orthogonal projection from ${\cal X}$ to
$\mathcal{H}$. The linearized Hamiltonian problem
(\ref{coupled-problem}) for non-zero eigenvalues $\lambda \neq 0$ is
rewritten as the generalized eigenvalue problem
\begin{equation}
\label{generalized_eigenvalue_problem} A u = \gamma K u, \qquad u\in
\mathcal{H},
\end{equation}
where $A = {\cal P} L_+ {\cal P}$, $K = {\cal P} L_-^{-1} {\cal P}$,
and $\gamma = - \lambda ^2$. We note that $K$ is a bounded
invertible self-adjoint operator on $\mathcal{H}$, while properties
of $A$ follow from those of $L_+$. Finitely many isolated
eigenvalues of the operators $A$ and $K^{-1}$ in ${\cal H}$ are
distributed between negative, zero and positive eigenvalues away of
$\sigma_c(L_{\pm})$. By the spectral theory of the self-adjoint
operators, the Hilbert space ${\cal H}$ can be equivalently
decomposed into two orthogonal sums of subspaces which are invariant
with respect to the operators $K$ and $A$:
\begin{eqnarray}
\label{decomposition-K} \mathcal{H} & = & \mathcal{H}_{K}^- \oplus
\mathcal{H}_{K}^+ \oplus \mathcal{H}_{K}^{\sigma_e(K)}, \\
\label{decomposition-A} \mathcal{H} & = & \mathcal{H}_{A}^- \oplus
\mathcal{H}_{A}^0 \oplus \mathcal{H}_{A}^+ \oplus
\mathcal{H}_{A}^{\sigma_e(A)},
\end{eqnarray}
where notation $-(+)$ stands for the negative (positive) isolated
eigenvalues, $0$ for the isolated kernel, and $\sigma_e$ for the
essential spectrum that includes the absolute continuous part and
embedded eigenvalues. It is clear that $\sigma_e(K)$ belongs to the
interval $(0,\omega_-^{-1}]$ and $\sigma_e(A)$ belongs to the
interval $[\omega_+,\infty)$. Since ${\cal P}$ is a projection
defined by eigenspaces of $L_-$ while $K = {\cal P} L_-^{-1} {\cal 
P}$, it is clear that ${\rm dim}(\mathcal{H}_{K}^-) = n(L_-)$ and 
${\rm dim}(\mathcal{H}_{K}^+) = p(L_-)$.

\begin{Proposition}
\label{lemma-constrained-space} Let $\omega_+ > 0$. There exist $n_0 
\geq 0$, $z_0 \geq 0$, and $z_1 \geq 0$, such that
\begin{equation}
\label{equality-constrained-space} {\rm dim}(\mathcal{H}_{A}^-) =
n(L_+) - n_0, \quad {\rm dim}(\mathcal{H}_{A}^0) = z(L_+) - z_0 +
z_1, \quad {\rm dim}(\mathcal{H}_{A}^+) \leq p(L_+) + n_0 + z_0 -
z_1.
\end{equation}
\end{Proposition}

\begin{Proof}
Let ${\rm ker}(L_-) = \{ v_1,v_2,...,v_n\} \in {\cal X}$ and define
the matrix-valued function $A(\mu)$:
$$
A_{ij}(\mu) = ((\mu-L_+)^{-1} v_i, v_j), \qquad 1 \leq i,j \leq n
$$
for all $\mu$ not in the spectrum of $L_+$. When $z(L_+) = 0$, the
first two equalities (\ref{equality-constrained-space}) follow by
the abstract Lemma 3.4 in \cite{CPV}, where $n_0$ is the number of
nonnegative eigenvalues of $A(0)$, $z_0 = 0$, and $z_1$ is the
number of zero eigenvalues of $A(0)$. When $z(L_+) \neq 0$, the same 
proof is extended for $A_0 = \lim\limits_{\mu \to 0^-} A(\mu)$ where 
$z_0$ is the number of eigenvectors in the kernel of $L_+$ in ${\cal 
X}$ which do not belong to the space ${\cal H}$ (see proof of 
Theorem 2.9 in \cite{CPV}). The last inequality in 
(\ref{equality-constrained-space}) is obtained by extending the 
analysis of \cite{CPV} from $\mu = 0$ to $\mu = \omega_+ > 0$, where 
the upper bound is achieved if all eigenvalues of $A_+ = 
\lim\limits_{\mu \to \omega_+^-} A(\mu)$ are negative.
\end{Proof}

Since $A$ has finitely many negative eigenvalues and $K$ has no
kernel in ${\cal H}$, there exists a small number $\delta > 0$ in
the gap $0 < \delta < |\sigma_{-1}|$, where $\sigma_{-1}$ is the
smallest (in absolute value) negative eigenvalue of $K^{-1} A$. The
operator $A + \delta K$ is continuously invertible in ${\cal H}$ and
the generalized eigenvalue problem
(\ref{generalized_eigenvalue_problem}) is rewritten in the shifted
form,
\begin{equation}
\label{shifted_eigenvalue_problem} (A + \delta K) u = (\gamma +
\delta)K u, \qquad u \in {\cal H}.
\end{equation}
By the spectral theory, an alternative decomposition of the Hilbert
space ${\cal H}$ exists for $0 < \delta < |\sigma_{-1}|$:
\begin{equation}
\label{decomposition-A-K-delta} \mathcal{H} = {\cal H}_{A + \delta
K}^- \oplus {\cal H}_{A + \delta K}^+ \oplus {\cal H}_{A+\delta
K}^{\sigma_e(A + \delta K)},
\end{equation}
where $\sigma_e(A + \delta K)$ belongs to the interval $[\omega_{A +
\delta K}, \infty)$ and $\omega_{A + \delta K}$ is the minimum of
$\sigma_c(A + \delta K)$. We will assume that $\omega_{A + \delta K}
> 0$ for $\delta > 0$ even if $\omega_+ = 0$.

\begin{Proposition}
\label{proposition-zero-splitting} Let the number of negative
(positive) eigenvalues of $A + \delta K$ bifurcating from the zero
eigenvalues of $A$ as $\delta > 0$ be denoted as $n_-$ ($n_+$).
When $\omega_+ > 0$, the splitting is complete so that
\begin{equation}
\label{splitting-zero-eigenvalue} {\rm dim}({\cal H}_{A + \delta
K}^{\pm}) = {\rm dim}({\cal H}_A^{\pm}) + n_{\pm}, \qquad {\rm
dim}({\cal H}_A^0) = n_- + n_+.
\end{equation}
When $\omega_+ = 0$ and ${\rm dim}({\cal H}_A^0) = 1$, the
statement remains valid provided that the bifurcating eigenvalue
of $A + \delta K$ is smaller than $\omega_{A + \delta K} > 0$.
\end{Proposition}

\begin{Proof}
For any $\delta \in (0,|\sigma_{-1}|)$, the operator $A + \delta K$ 
has no zero eigenvalues in ${\cal H}$. The equality 
(\ref{splitting-zero-eigenvalue}) follows from the definition of 
$n_-$ and $n_+$.
\end{Proof}

The Pontryagin's Invariant Subspace theorem can be applied to the
product of two bounded invertible self-adjoint operators $(A +
\delta K)^{-1}$ and $K$ in ${\cal H}$.

\section{The proof of the Pontryagin's Invariant Subspace theorem}

We shall develop an abstract theory of Pontryagin spaces with
sign-indefinite metric, where the main result is the Pontryagin's
Invariant Subspace theorem.

\begin{Definition}
Let $\mathcal{H}$ be a Hilbert space equipped with the inner
product $(\cdot,\cdot)$ and the sesquilinear form
$[\cdot,\cdot]$\footnote{We say that a complex-valued form $[u,v]$
on the product space $\mathcal{H} \times \mathcal{H}$ is a
sesquilinear form if it is linear in $u$ for each fixed $v$ and
linear with complex conjugate in $v$ for each fixed $u$.}. The
Hilbert space $\mathcal{H}$ is called the Pontryagin space
(denoted as $\Pi_{\kappa}$) if it can be decomposed into the sum,
which is orthogonal with respect to $[\cdot,\cdot]$,
\begin{equation} \label{decomposition-Pi}
\mathcal{H} \doteq \Pi_{\kappa} = \Pi_{+} \oplus \Pi_{-},
\end{equation}
where $\Pi_{+}$ is a Hilbert space with the inner product
$(\cdot,\cdot) = [\cdot,\cdot]$, $\Pi_{-}$ is a Hilbert space with
the inner product $(\cdot,\cdot) = -[\cdot,\cdot]$, and $\kappa =
{\rm dim}(\Pi_{-}) < \infty$.
\end{Definition}

\begin{Remark}
{\rm We shall write components of an element $x$ in the Pontryagin
space $\Pi_{\kappa}$ as a vector $ x = \{x_-,x_+\}$. The
orthogonal sum (\ref{decomposition-Pi}) implies that any non-zero
element $x \neq 0$ is represented by two terms,
\begin{equation}
\label{decomposition-Pontryagin} \forall x \in \Pi_{\kappa} : \quad
x = x_+ + x_-,
\end{equation}
such that
\begin{equation}
\label{properties-decomposition} [x_+,x_-] = 0, \quad [x_+,x_+]
> 0, \quad [x_-,x_-] < 0,
\end{equation}
and $\Pi_+ \cap \Pi_- = \varnothing$.}
\end{Remark}

\begin{Definition}
We say that $\Pi$ is a non-positive subspace of $\Pi_{\kappa}$ if
$[x,x] \leq 0$ $\forall x \in \Pi$. We say that the non-positive
subspace $\Pi$ has the maximal dimension $\kappa$ if any subspace of
$\Pi_{\kappa}$ of dimension higher than $\kappa$ is not a 
non-positive subspace of $\Pi_{\kappa}$. Similarly, $\Pi$ is a 
non-negative (neutral) subspace of $\Pi_{\kappa}$ if $[x,x] \geq 0$ 
($[x,x] = 0$) $\forall x \in \Pi$. The sign of $[x,x]$ on the 
element $x$ of a subspace $\Pi$ is called the Krein signature of the 
subspace $\Pi$.
\end{Definition}

\begin{Theorem}[Pontryagin]
\label{theorem-Pontryagin} Let $T$ be a self-adjoint bounded
operator in $\Pi_{\kappa}$, such that $[T \cdot, \cdot] = [\cdot, T
\cdot]$. There exists a $T$-invariant non-positive subspace of
$\Pi_{\kappa}$ of the maximal dimension $\kappa$.
\end{Theorem}

\begin{Remark}
{\rm There are historically two completely different approaches to
the proof of this theorem. A proof based on theory of analytic
functions was given by Pontryagin \cite{Pontryagin} while a proof
based on angular operators was given by Krein \cite{krein} and later
developed by students of M.G. Krein \cite{Azizov,iohvidov}. Theorem
\ref{theorem-Pontryagin} was rediscovered by Grillakis \cite{G90}
with the use of topology. We will describe a geometric proof of
Theorem \ref{theorem-Pontryagin} based on the Fixed Point theorem.
The proof uses the Cayley transformation of a self-adjoint operator
in $\Pi_{\kappa}$ to a unitary operator in $\Pi_{\kappa}$ (Lemma
\ref{invariant}) and the Krein's representation of the maximal
non-positive subspace of $\Pi_{\kappa}$ in terms of a graph of the
contraction map (Lemma \ref{contraction}). While many statements of
our analysis are available in the literature, details of the proofs
are missing. Our presentation gives full details of the proof of
Theorem \ref{theorem-Pontryagin} (see \cite{GKP} for a similar
treatment in the case of compact operators). }
\end{Remark}

\begin{Lemma}
\label{invariant} Let $T$ be a linear operator in $\Pi_{\kappa}$ and
$z \in \mathbb{C}$, $\Im(z) > 0$ be a regular point of the operator
$T$, such that $z \in \rho (T)$. Let $U$ be the Cayley transform of
$T$ defined by $U = (T - \bar{z})(T - z)^{-1}$. The operators $T$
and $U$ have the same invariant subspaces in $\Pi_{\kappa}$.
\end{Lemma}

\begin{Proof}
Let $\Pi$ be a finite-dimensional invariant subspace of the operator
$T$ in $\Pi_{\kappa}$. It follows from $z \in \rho(T)$ that $(T - z)
\Pi = \Pi$ then $(T - z)^{-1} \Pi = \Pi$ and $(T - \bar{z)}(T -
z)^{-1} \Pi \subseteq \Pi$, i.e $U \Pi \subseteq \Pi$. Conversely,
let $\Pi$ be an invariant subspace of the operator $U$. It follows
from $U - I = (z - \bar{z})(T - z)^{-1}$ that $1 \in \rho(U)$
therefore $\Pi = (U - I)\Pi = (T - z)^{-1} \Pi$. From there, $\Pi
\subseteq {\rm dom}(T)$ and $(T - z) \Pi = \Pi$ so $T \Pi \subseteq
\Pi$.
\end{Proof}

\begin{Corollary}
\label{corollary-TU} If $T$ is a self-adjoint operator in
$\Pi_{\kappa}$, then $U$ is a unitary operator in $\Pi_{\kappa}$.
\end{Corollary}

\begin{Proof}
We shall prove that $[U g, Ug] = [g,g]$, where $g \in {\rm dom}(U)$,
by the explicit computation:
\begin{eqnarray*}
[U g, U g] = \left[(T - \bar{z}) f, (T - \bar{z}) f
\right] & = & [Tf,Tf] - z [f,Tf] - \bar{z} [T f, f]  + |z|^2[f,f], \\
\left[ g, g \right] = \left[ (T - z) f, (T - z) f \right] & = &
[Tf,Tf] - \bar{z}[f,Tf] - z [Tf,f] + |z|^2 [f,f],
\end{eqnarray*}
where we have introduced $f \in {\rm dom}(T)$ such that $f =
(T-z)^{-1} g$.
\end{Proof}

\begin{Lemma}
\label{contraction} A linear subspace $\Pi \subseteq \Pi_{\kappa}$
is a $\kappa$-dimensional non-positive subspace of $\Pi_{\kappa}$
if and only if it is a graph of the contraction map ${\cal K} :
\Pi_{-} \rightarrow \Pi_{+}$, such that $\Pi = \{ x_-,{\cal K} x_-
\}$ and $\| {\cal K} x_- \| \leq \| x_- \|$.
\end{Lemma}

\begin{Proof}
Let $\Pi = \{x_-,x_+\}$ be a $\kappa$-dimensional non-positive
subspace of $\Pi_{\kappa}$. We will show that there exist a
contraction map ${\cal K} : \Pi_- \mapsto \Pi_+$ such that $\Pi$ is 
a graph of ${\cal K}$. Indeed, the subspace $\Pi$ is a graph of a 
linear operator ${\cal K}$ if and only if it follows from $\{0,x_+\} 
\in \Pi$ that $x_+ = 0$. Since $\Pi$ is non-positive with respect to 
$[\cdot,\cdot]$, then $[x,x] = \| x_+ \|^2 - \| x_- \|^2 \leq 0$, 
where $\| \cdot \|$ is a norm in $\cal H$. As a result, $0 \leq \| 
x_+ \| \leq \| x_- \|$ and if $x_- = 0$ then $x_+ = 0$. Moreover, 
for any $x_- \in \Pi_-$, it is true that $\| {\cal K} x_- \| \leq \| 
x_- \|$ such that ${\cal K}$ is a contraction map. Conversely, let 
${\cal K}$ be a contraction map ${\cal K} : \Pi_- \mapsto \Pi_+$. 
The graph of ${\cal K}$ belongs to the non-positive subspace of 
$\Pi_{\kappa}$ as
$$
[x,x] = \| x_+ \|^2 - \| x_- \|^2 = \| {\cal K} x_- \|^2 - \| x_-
\|^2 \leq 0.
$$
Let $\Pi = \{ x_-, {\cal K} x_- \}$. Since ${\rm dim}(\Pi_-) =
\kappa$, then ${\rm dim}(\Pi) = \kappa$. \footnote{Extending
arguments of Lemma \ref{contraction}, one can prove that the
subspace $\Pi$ is strictly negative with respect to
$[\cdot,\cdot]$ if and only if it is a graph of the strictly
contraction map ${\cal K} : \Pi_- \mapsto \Pi_+$, such that $\Pi =
\{ x_-,{\cal K} x_- \}$ and $\| {\cal K} x_- \| < \| x_- \|$.}
\end{Proof}

\begin{Proof1}{\bf of Theorem \ref{theorem-Pontryagin}}.
Let $z \in \C$ and $\Im(z) > 0$. Then, $z$ is a regular point of the
self-adjoint operator $T$ in $\Pi_{\kappa}$. Let $U =(T - \bar{z})(T
- z)^{-1}$ be the Cayley transform of $T$. By Corollary
\ref{corollary-TU}, $U$ is a unitary operator in $\Pi_{\kappa}$. By
Lemma \ref{invariant}, $T$ and $U$ have the same invariant subspaces
in $\Pi_{\kappa}$. Therefore, the existence of the maximal
non-positive invariant subspace for the self-adjoint operator $T$
can be proved from the existence of such a subspace for the unitary
operator $U$. Let $x = \{ x_-,x_+ \}$ and
$$
U = \left[ \begin{array}{cc} U_{11} & U_{12}\\ U_{21} & U_{22}
\end{array} \right]
$$
be the matrix representation of the operator $U$ with respect to
the decomposition (\ref{decomposition-Pi}). Let $\Pi$ denote a
$\kappa$-dimensional non-positive subspace in $\Pi_{\kappa}$.
Since $U$ has an empty kernel in $\Pi_{\kappa}$ and $U$ is unitary
in $\Pi_{\kappa}$ such that $[U x_-, U x_-] = [x_-,x_-] \leq 0$,
then $\tilde{\Pi} = U \Pi$ is also a $\kappa$-dimensional
non-positive subspace of $\Pi_{\kappa}$. By Lemma
\ref{contraction}, there exist two contraction mappings ${\cal K}$
and $\tilde{\cal K}$ for subspaces $\Pi$ and $\tilde{\Pi}$,
respectively. Therefore, the assignment $\tilde{\Pi} = U \Pi$ is
equivalent to the system,
$$
\left( \begin{array}{c} \tilde{x}_- \\ \tilde{\cal K} \tilde{x}_-
\end{array}  \right) = \left[ \begin{array}{cc} U_{11} & U_{12}\\
U_{21} & U_{22} \end{array} \right] \left( \begin{array}{c} x_- \\
{\cal K} x_- \end{array}  \right) = \left(
\begin{array}{c}
(U_{11} + U_{12}{\cal K}) x_-\\ (U_{21} + U_{22} {\cal K})x_-
\end{array} \right),
$$
such that
$$
U_{21} + U_{22} {\cal K} = \tilde{\cal K}(U_{11} + U_{12} {\cal K}).
$$
We shall prove that the operator $(U_{11} + U_{12} {\cal K})$ is
invertible. By contradiction, we assume that there exists $x_-
\neq 0$ such that $\tilde{x}_- = (U_{11} + U_{12} {\cal K}) x_- =
0$. Since $\tilde{x}_- = 0$ implies that $\tilde{x}_+ =
\tilde{\cal K} \tilde{x}_- = 0$, we obtain that $[x_-,{\cal K}
x_-]^T$ is an eigenvector in the kernel of $U$. However, $U$ has
an empty kernel in $\Pi_{\kappa}$ such that $x_- = 0$. Let
$F({\cal K})$ be an operator-valued function in the form,
$$
F({\cal K}) = (U_{21} + U_{22} {\cal K})(U_{11} + U_{12} {\cal
K})^{-1},
$$
such that $\tilde{\cal K} = F({\cal K})$. By Lemma
\ref{contraction}, the operator $F({\cal K})$ maps the operator
unit ball $\| {\cal K} \| \leq 1$ to itself. Since $U$ is a
continuous operator and $U_{12}$ is a finite-dimensional operator,
then $U_{12}$ is a compact operator. Hence the operator ball $\|
{\cal K} \| \leq 1$ is a weakly compact set and the function
$F({\cal K})$ is continuous with respect to weak topology. By the
Schauder's Fixed-Point Principle, there exists a fixed point
${\cal K}_0$ such that $F({\cal K}_0) = {\cal K}_0$ and $\|{\cal
K}_0 \| \leq 1$. By Lemma \ref{contraction}, the graph of ${\cal
K}_0$ defines the $\kappa$-dimensional non-positive subspace
$\Pi$, which is invariant with respect to $U$.

It remains to prove that the $\kappa$-dimensional non-positive
subspace $\Pi$ has the maximal dimension that a non-positive 
subspace of $\Pi_{\kappa}$ can have. By contradiction, we assume 
that there exists a $(\kappa + 1)$-dimensional non-positive subspace 
$\tilde{\Pi}$. Let $\{e_1,e_2,...,e_{\kappa} \}$ be a basis in 
$\Pi_-$ in the canonical decomposition 
(\ref{decomposition-Pontryagin}). We fix two elements $y_1, y_2 \in 
\tilde{\Pi}$ with the same projections to 
$\{e_1,e_2,...,e_{\kappa}\}$, such that
\begin{eqnarray*}
y_1 & = & \alpha_1 e_1 + \alpha_2 e_2 + ... + \alpha_{\kappa}
e_{\kappa} + y_{1p}, \\
y_2 & = & \alpha_1 e_1 + \alpha_2 e_2 + ... + \alpha_{\kappa}
e_{\kappa} + y_{2p},
\end{eqnarray*}
where $y_{1p}, y_{2p} \in \Pi_{+}$. It is clear that $y_1 -y_2 =
y_{1p} - y_{2p} \in \Pi_+$ such that $[y_{1p}-y_{2p},y_{1p}-y_{2p}] 
> 0$. On the other hand, $y_1 - y_2 \in \tilde{\Pi}$, such that $[
y_1-y_2,y_1-y_2] \leq 0$. Hence $y_{1p} = y_{2p}$ and then $y_1 
=y_2$. We proved that any vector in $\tilde{\Pi}$ is uniquely 
determined by the $\kappa$-dimensional basis 
$\{e_1,e_2,...,e_{\kappa}\}$ and the non-positive subspace 
$\tilde{\Pi}$ is hence $\kappa$-dimensional.
\end{Proof1}

\section{Bounds on eigenvalues of the generalized eigenvalue problem}

We shall apply Theorem \ref{theorem-Pontryagin} to the product of
two bounded invertible self-adjoint operators $B = (A + \delta
K)^{-1}$ and $K$, where $\delta \in (0,|\sigma_{-1}|)$ and
$\sigma_{-1}$ is the smallest negative eigenvalue of $K^{-1}A$.
Properties of self-adjoint operators $A$ and $K^{-1}$ in ${\cal H}$ 
follow from properties of self-adjoint operators $L_{\pm}$ in ${\cal 
X}$, which are summarized in the main assumptions P1--P2. With a 
slight abuse of notations, we shall denote eigenvalues of the 
operator $T = B K$ by $\lambda$, which is expressed in terms of the 
eigenvalue $\gamma$ of the shifted generalized eigenvalue problem 
(\ref{shifted_eigenvalue_problem}) by $\lambda = (\gamma + 
\delta)^{-1}$. We note that $\lambda$ here does not correspond to 
$\lambda$ used in the linearized Hamiltonian problem 
(\ref{coupled-problem}).

\begin{Lemma}
\label{BG product} Let $\mathcal{H}$ be a Hilbert space with the
inner product $(\cdot,\cdot)$ and $B,K : \mathcal{H} \rightarrow
\mathcal{H}$ be bounded invertible self-adjoint operators in
$\mathcal{H}$. Define the sesquilinear form $[\cdot, \cdot] =
(K\cdot,\cdot)$ and extend $\mathcal{H}$ to the Pontryagin space
$\Pi_{\kappa}$, where $\kappa$  is the finite number of negative
eigenvalues of $K$ counted with their multiplicities. The operator
$T = B K$ is self-adjoint in $\Pi_{\kappa}$ and there exists a
$\kappa$-dimensional non-positive subspace of $\Pi_{\kappa}$ which
is invariant with respect to $T$.
\end{Lemma}

\begin{Proof}
It follows from the orthogonal sum decomposition
(\ref{decomposition-K}) that the quadratic form $(K \cdot,\cdot)$ is
strictly negative on the $\kappa$-dimensional subspace ${\cal
H}_K^-$ and strictly positive on the infinite-dimensional subspace
${\cal H}_K^+ \oplus {\cal H}_K^{\sigma_e(K)}$. By continuity and
Gram--Schmidt orthogonalization, the Hilbert space $\mathcal{H}$ is
extended to the Pontryagin space $\Pi_{\kappa}$ with respect to the
sesquilinear form $[\cdot, \cdot] = (K\cdot,\cdot)$. The bounded
operator $T = B K$ is self-adjoint in $\Pi_{\kappa}$, since $B$ and
$K$ are self-adjoint in ${\cal H}$ and
$$
[T\cdot,\cdot] = (K B K \cdot,\cdot) = (K \cdot, B K \cdot) =
[\cdot,T\cdot].
$$
Existence of the $\kappa$-dimensional non-positive $T$-invariant
subspace of $\Pi_{\kappa}$ follows from Theorem
\ref{theorem-Pontryagin}.
\end{Proof}

\begin{Remark}
{\rm The decomposition (\ref{decomposition-Pi}) of the Pontryagin
space $\Pi_{\kappa}$ is canonical in the sense that $\Pi_+ \cap
\Pi_- = \varnothing$.  We will now consider various sign-definite
subspaces of $\Pi_{\kappa}$ which are invariant with respect to
the self-adjoint operator $T = B K$ in $\Pi_{\kappa}$. In general,
these invariant sign-definite subspaces do not provide a canonical
decomposition of $\Pi_{\kappa}$. Let us denote the invariant
subspace of $T$ associated with complex eigenvalues in the upper
(lower) half-plane as $\mathcal{H}_{c^+}$ ($\mathcal{H}_{c^-}$)
and the non-positive (non-negative) invariant subspace of $T$
associated with real eigenvalues as $\mathcal{H}_n
(\mathcal{H}_p)$. The invariant subspace $\mathcal{H}_n$,
prescribed by Lemma \ref{BG product},  may include both isolated
and embedded eigenvalues of $T$ in $\Pi_{\kappa}$. We will show
that this subspace does not include the residual and absolutely
continuous parts of the spectrum of $T$ in $\Pi_{\kappa}$. }
\end{Remark}

\subsection{Residual and absolutely continuous spectra of $T$ in
$\Pi_{\kappa}$}

\begin{Definition}
We say that $\lambda$ is a point of the residual spectrum of $T$ in 
$\Pi_{\kappa}$ if ${\rm Ker}(T-\lambda I) = \varnothing$ but 
$\overline{{\rm Range}(T - \lambda I)} \neq \Pi_{\kappa}$ and 
$\lambda$ is a point of the continuous spectrum of $T$ in 
$\Pi_{\kappa}$ if ${\rm Ker}(T-\lambda I) = \varnothing$ but ${\rm 
Range}(T - \lambda I) \neq \overline{{\rm Range}(T - \lambda I)} =  
\Pi_{\kappa}$.
\end{Definition}

\begin{Lemma}
\label{lemma-residual} No residual part of the spectrum of $T$ in
$\Pi_{\kappa}$ exists.
\end{Lemma}

\begin{Proof}
By a contradiction, assume that $\lambda$ belongs to the residual
part of the spectrum of $T$ in $\Pi_{\kappa}$ such that ${\rm
Ker}(T - \lambda I) =  \varnothing$ but ${\rm Range}(T - \lambda
I)$ is not dense in $\Pi_{\kappa}$. Let $g \in \Pi_{\kappa}$ be
orthogonal to ${\rm Range}(T - \lambda I)$, such that
$$
\forall f \in \Pi_{\kappa} : \quad 0 = [(T - \lambda I)f, g] = [
f, (T - \bar{\lambda} I) g].
$$
Therefore, $(T - \bar{\lambda} I) g = 0$, that is $\bar{\lambda}$
is an eigenvalue of $T$. By symmetry of eigenvalues, $\lambda$ is
also an eigenvalue of $T$ and hence it can not be in the residual
part of the spectrum of $T$.
\end{Proof}

\begin{Remark}
{\rm It is assumed in \cite{G90} that the residual part of
spectrum is empty and that the kernels of operators $A$ and $K$
are empty. The first assumption is now proved in Lemma
\ref{lemma-residual} and the second assumption is removed with the
use of the shifted generalized eigenvalue problem
(\ref{shifted_eigenvalue_problem}).  }
\end{Remark}

\begin{Lemma}
The absolutely continuous part of the spectrum of $T$ in
$\Pi_{\kappa}$ is real.
\end{Lemma}

\begin{Proof}
Let $P^{+}$ and $P^{-}$ be orthogonal projectors to $\Pi^{+}$ and
$\Pi^{-}$ respectively, such that $I = P^+ + P^-$. Since
$\Pi^{\pm}$ are defined by $(\cdot, K \cdot)$, the self-adjoint
operator $K$ admits the polar decomposition $K = J|K|$, where $J =
P^{+} - P^{-}$ and $|K|$ is a positive operator. The operator $T =
B K$ is similar to the operator
$$
|K|^{1/2}B J|K|^{1/2} = |K|^{1/2} B J|K|^{1/2}( J + 2 P^{-} ) =
|K|^{1/2} B |K|^{1/2} + 2 |K|^{1/2} B J|K|^{1/2} P^{-}.
$$
Since $P^-$ is a projection to a finite-dimension subspace, the
operator $|K|^{1/2} B J|K|^{1/2}$ is the finite-dimensional
perturbation of the self-adjoint operator $|K|^{1/2} B |K|^{1/2}$.
By perturbation theory \cite{Kato}, the absolutely continuous part
of the spectrum of the self-adjoint operator $|K|^{1/2} B
|K|^{1/2}$ is the same as that of $|K|^{1/2} B J|K|^{1/2}$. By
similarity transformation, it is the same as that of $T$.
\end{Proof}

\begin{Lemma}[Cauchy-Schwartz]
\label{cauchy-inequality} Let $\Pi$ be either non-positive or
non-negative subspace of $\Pi_{\kappa}$.  Then,
\begin{equation}
\label{cauchy-schwartz} \forall f,g \in \Pi : \quad |[f,g]|^2 \leq
[f,f][g,g].
\end{equation}
\end{Lemma}

\begin{Proof}
The proof resembles that of the standard Cauchy--Schwartz
inequality. Let $\Pi$ be a non-positive subspace of
$\Pi_{\kappa}$, Then, for any $f,g \in \Pi$ and any $\alpha,\beta
\in \mathbb{C}$, we have $[\alpha f + \beta g,\alpha f + \beta g]
\leq 0$ and
\begin{equation}
\label{help-CS} [\alpha f + \beta g,\alpha f + \beta g]  = [f,f]
|\alpha|^2 + [f,g] \alpha \bar{\beta} + [g,f]\bar{\alpha} \beta
+[g,g]|\beta|^2.
\end{equation}
If $[f,g] = 0$, then (\ref{cauchy-schwartz}) is satisfied as $[f,f] 
\leq 0$ and $[g,g] \leq 0$. If $[f,g] \neq 0$, let $\alpha \in 
\mathbb{R}$ and 
$$
\beta = \frac{[f,g]}{|[f,g]|},
$$
such that 
$$
[f,f] \alpha^2 + 2 \alpha |[f,g]| + [g,g] \leq 0.
$$
The last condition is satisfied if the discriminant of the quadratic 
equation is non-positive such that $4 |[f,g]|^2 - 4 [f,f][g,g] \leq 
0$, that is (\ref{cauchy-schwartz}). Let $\Pi$ be a non-negative 
subspace of $\Pi_{\kappa}$.  Then, for any $f,g \in \Pi$ and any 
$\alpha,\beta \in \mathbb{C}$, we have $[\alpha f + \beta g,\alpha f 
+ \beta g] \geq 0$ and the same arguments result in the same 
inequality (\ref{cauchy-schwartz}).
\end{Proof}

\begin{Corollary}
\label{corollary-cauchy} Let $\Pi$ be either non-positive or
non-negative subspace of $\Pi_{\kappa}$. Let $f \in \Pi$ such that
$[f,f] =0$. Then $[f,g] = 0$ $\forall g \in \Pi$.
\end{Corollary}

\begin{Proof}
The proof follows from (\ref{cauchy-schwartz}) since $0 \leq
|[f,g]|^2 \leq 0$.
\end{Proof}

\begin{Lemma}
\label{lemma-cauchy-invariance} Let $\Pi$ be an invariant subspace
of $\Pi_{\kappa}$ with respect to operator $T$ and $\Pi^{\perp}$
be the orthogonal compliment of $\Pi$ in $\Pi_{\kappa}$ with
respect to $[\cdot,\cdot]$. Then, $\Pi^{\perp}$ is also invariant
with respect to $T$.
\end{Lemma}

\begin{Proof}
For all $f \in {\rm Dom}(T) \cap \Pi$, we have  $T f \in \Pi$. Let
$g \in {\rm Dom}(T) \cap \Pi^{\perp}$. Then $[g,Tf] = [Tg,f] = 0$.
\end{Proof}

\begin{Theorem}
Let $\Pi_c$ be a subspace related to the absolute continuous
spectrum of $T$ in $\Pi_{\kappa}$. Then $[f,f] > 0$ $\forall f \in
\Pi_c$.
\end{Theorem}

\begin{Proof}
By contradiction, assume that there exists $f_0 \in \Pi_c$ such that 
$[f_0,f_0] < 0$. Since $\Pi_c$ is a subspace for the absolutely 
continuous spectrum, there exists a continuous family of functions 
$f_{\alpha} \in \Pi_c$ such that $[f_{\alpha},f_{\alpha}] < 0$ and 
hence $f_{\alpha} \in \Pi_-$ in the decomposition 
(\ref{decomposition-Pi}). However, this contradicts to the fact that 
${\rm dim}(\Pi_-) = \kappa < \infty$. Therefore, $\Pi_c$ is a 
non-negative subspace of $\Pi_{\kappa}$. Assume that there exists an 
element $f_0 \in \Pi_c$ such that $[f_0,f_0] = 0$. By Corollary 
\ref{corollary-cauchy}, $f_0 \in \Pi_c^{\perp}$. By Lemma 
\ref{lemma-cauchy-invariance}, $\Pi_c^{\perp}$ is an invariant 
subspace of $\Pi_{\kappa}$. The intersection of invariant subspaces 
is invariant, such that $\Pi_c \cap \Pi_c^{\perp}$ is a neutral 
invariant subspace of $\Pi_{\kappa}$. By Lemma \ref{BG product}, 
${\rm dim}(\Pi_c \cap \Pi_c^{\perp}) \leq \kappa$. Therefore, $f_0$ 
is an element of a finite-dimensional invariant subspace of 
$\Pi_{\kappa}$, which is a contradiction to the fact that $f_0 \in
\Pi_c$.
\end{Proof}

\subsection{Isolated and embedded eigenvalues of $T$ in
$\Pi_{\kappa}$}

\begin{Definition}
\label{definition-multiple} We say that $\lambda$ is an eigenvalue
of $T$ in $\Pi_{\kappa}$ if ${\rm Ker}(T - \lambda I) \neq
\varnothing$. The eigenvalue of $T$ is said to be semi-simple if it 
is not simple and the algebraic and geometric multiplicities 
coincide. Otherwise, the non-simple eigenvalue is said to be 
multiple. Let $\lambda_0$ be an eigenvalue of $T$ with algebraic 
multiplicity $n$ and geometric multiplicity one. The canonical basis 
for the corresponding $n$-dimensional eigenspace of $T$ is defined 
by the Jordan chain of eigenvectors,
\begin{equation}
\label{Jordan-chain-formula} f_j \in \Pi_{\kappa} : \quad T f_j =
\lambda_0 f_j + f_{j-1}, \qquad j = 1,...,n,
\end{equation}
where $f_0 = 0$.
\end{Definition}

\begin{Lemma}
\label{ortspaces} Let $\mathcal H_{\lambda}$ and $\mathcal
H_{\mu}$ be eigenspaces of the eigenvalues $\lambda$ and $\mu$ of
the operator $T$ in $\Pi_{\kappa}$ and $\lambda \neq \bar{\mu}$.
Then $\mathcal H_{\lambda}$ is orthogonal  to $\mathcal H_{\mu}$
with respect to $[\cdot,\cdot]$.
\end{Lemma}

\begin{Proof}
Let $n$ and $m$ be dimensions of ${\cal H}_{\lambda}$ and ${\cal
H}_{\mu}$, respectively, such that $n \geq 1$ and $m \geq 1$. By
Definition \ref{definition-multiple}, it is clear that
\begin{eqnarray}
\label{help-system-1} f \in {\cal H}_{\lambda} \Longleftrightarrow
(T - \lambda I)^n f = 0, \\ \label{help-system-2} g \in {\cal
H}_{\mu} \Longleftrightarrow (T - \mu I)^m g = 0.
\end{eqnarray}
We should prove that $[f,g] = 0$ by induction for $n+m \geq 2$. If
$n+m = 2$ ($n = m = 1$), then it follows from the system
(\ref{help-system-1})--(\ref{help-system-2}) that
$$
(\lambda - \bar{\mu}) [f,g] = 0, \qquad f \in {\cal H}_{\lambda},
\quad g \in {\cal H}_{\mu},
$$
such that $[f,g] = 0$ for $\lambda \neq \bar{\mu}$. Let us assume
that subspaces $\mathcal H_{\lambda}$ and $\mathcal H_{\mu}$ are
orthogonal for $2 \leq n + m \leq k$ and prove that extended
subspaces $\tilde{\mathcal H}_{\lambda}$ and $\tilde{\mathcal
H}_{\mu}$ remain orthogonal for $\tilde{n} = n+1$, $\tilde{m} = m$
and $\tilde{n} = n$, $\tilde{m} = m + 1$. In either case, we define
$\tilde{f} =(T - \lambda I) f$ and $\tilde{g} =(T - \mu I)g$, such
that
\begin{eqnarray*}
f \in \tilde{\cal H}_{\lambda} \Longleftrightarrow (T - \lambda
I)^{\tilde{n}} f = (T - \lambda I)^n \tilde{f} = 0, \\ g \in
\tilde{\cal H}_{\mu} \Longleftrightarrow (T - \mu I)^{\tilde{m}} g
= (T - \mu I)^m \tilde{g} = 0.
\end{eqnarray*}
By the inductive assumption, we have $[\tilde{f},g] = 0$ and
$[f,\tilde{g}] = 0$ in either case, such that
\begin{equation}
\label{ur3} [(T - \lambda I)f,g] = 0  \qquad [f,(T - \mu I)g] = 0.
\end{equation}
By using the system (\ref{help-system-1})--(\ref{help-system-2}) and
the relations (\ref{ur3}), we obtain that
$$
(\lambda - \bar{\mu})[f,g] = 0,  \qquad f \in \tilde{\cal
H}_{\lambda}, \quad g \in \tilde{\cal H}_{\mu},
$$
from which the statement follows by the induction method.
\end{Proof}

\begin{Lemma}
\label{Jordan chain} Let ${\mathcal H}_{\lambda_0}$ be eigenspace of
a multiple real isolated eigenvalue $\lambda_0$ of $T$ in
$\Pi_{\kappa}$ and $\{f_1,f_2,...f_n\}$ be the Jordan chain of
eigenvectors. Let ${\cal H}_0 = {\rm span}\{f_1,f_2,...,f_{k}\}
\subset {\cal H}_{\lambda_0}$, where $k = {\rm round}(n/2)$, and
$\tilde{\cal H}_0 = {\rm span}\{f_1,f_2,...,f_{k},f_{k+1} \} \subset
{\cal H}_{\lambda_0}$.
\begin{itemize}
\item If $n$ is even ($n = 2k$), the neutral subspace ${\cal H}_0$
is the maximal sign-definite subspace of ${\cal H}_{\lambda_0}$.

\item If $n$ is odd ($n = 2k+1$), the subspace $\tilde{\cal H}_0$
is the maximal non-negative subspace of $\Pi_{\kappa}$ if
$[f_1,f_n] > 0$ and the maximal non-positive subspace of
$\Pi_{\kappa}$ if $[f_1,f_n] < 0$, while the neutral subspace
${\cal H}_0$ is the maximal non-positive subspace of
$\Pi_{\kappa}$ if $[f_1,f_n] > 0$ and the maximal non-negative
subspace of $\Pi_{\kappa}$ if $[f_1,f_n] < 0$.
\end{itemize}
\end{Lemma}

\begin{Proof}
Without loss of generality we will consider the case $\lambda_0 = 0$
(if $\lambda_0 \neq 0$ the same argument is applied to the shifted
self-adjoint operator $\tilde{T} = T - \lambda_0 I$). We will show 
that $[f,f] = 0$ $\forall f \in {\cal H}_0$. By a decomposition over 
the basis in ${\cal H}_0$, we obtain
\begin{equation}
\label{decomposition-f-g} [f,f] = \sum_{i=1}^k \sum_{j=1}^k \alpha_i
\bar{\alpha}_j \left[  f_i,  f_j \right].
\end{equation}
We use that
$$
[f_i,f_j] = [ T f_{i+1},Tf_{j+1}] = ... = \left[ T^k f_{i+k}, T^k
f_{j+k} \right] = \left[T^{2k} f_{i+k}, f_{j+k} \right],
$$
for any $1 \leq i,j \leq k$. In case of even $n = 2k$, we have 
$[f_i,f_j] = [T^n f_{i+k},f_{j+k}] = 0$ for all $1 \leq i,j \leq k$. 
In case of odd $n = 2k + 1$, we have $[f_i,f_j] = [T^{n+1} 
f_{i+k+1},f_{j+k+1}] = 0$ for all $1 \leq i,j \leq k$. Therefore, 
${\cal H}_0$ is a neutral subspace of ${\cal H}_{\lambda_0}$. To 
show that it is actually the maximal neutral subspace of ${\cal 
H}_{\lambda_0}$, let ${\cal H}'_0 = {\rm span}\{f_1,f_2,...,f_{k}, 
f_{k_0}\}$, where $k+1 \leq k_0 \leq n$. Since $f_{n+1}$ does not 
exist in the Jordan chain (\ref{Jordan-chain-formula}) (otherwise, 
the algebraic multiplicity is $n+1$) and $\lambda_0$ is an isolated 
eigenvalue, then $[f_1,f_n] \neq 0$ by the Fredholm theory. It 
follows from the Jordan chain (\ref{Jordan-chain-formula}) that
\begin{equation}
\label{help-help} [f_1,f_n] = [T^{m-1} f_m, f_n] = [f_m, T^{m-1} 
f_n] = [f_m,f_{n-m+1}] \neq 0.
\end{equation}
When $n = 2k$, we have $1 \leq n-k_0+1 \leq k$, such that 
$[f_{k_0},f_{n-k_0+1}] \neq 0$ and the subspace ${\cal H}_0'$ is 
sign-indefinite in the decomposition (\ref{decomposition-f-g}). When 
$n = 2k+1$, we have $1 \leq n-k_0+1 \leq k$ for $k_0 \geq k+2$ and 
$n-k_0+1 = k+1$ for $k_0 = k+1$. In either case, 
$[f_{k_0},f_{n-k_0+1}] \neq 0$ and the subspace ${\cal H}_0'$ is 
sign-indefinite in the decomposition (\ref{decomposition-f-g}) 
unless $k_0 = k+1$. In the latter case, we have $[f_{k+1},f_{k+1}] = 
[f_1,f_n] \neq 0$ and $[f_j,f_{k+1}] = [T^{2k} f_{j+k},f_n] = 0$ for 
$1 \leq j \leq k$. As a result, the subspace $\tilde{\cal H}_0 
\equiv {\cal H}_0'$ is non-negative for $[f_1,f_n] 
> 0$ and non-positive for $[f_1,f_n] < 0$.
\end{Proof}

\begin{Corollary}
\label{positive(negative) combination} Let ${\mathcal
H}_{\lambda_0}$ be eigenspace of a semi-simple real isolated
eigenvalue $\lambda_0$ of $T$ in $\Pi_{\kappa}$. Then, the
eigenspace ${\mathcal H}_{\lambda_0}$ is either strictly positive
or strictly negative subspace of $\Pi_{\kappa}$ with respect to
$[\cdot,\cdot]$.
\end{Corollary}

\begin{Proof}
The proof follows by contradiction. Let $\tilde{f}$ be a particular
linear combination of eigenvectors in ${\cal H}_{\lambda_0}$, such
that $[\tilde{f},\tilde{f}] = 0$. By the Fredholm theory, there 
exists a solution of the Jordan chain (\ref{Jordan-chain-formula}), 
where $f_1 \equiv \tilde{f}$. Then, the eigenvalue $\lambda_0$ could 
not be semi-simple.
\end{Proof}

\begin{Remark}
\label{embedded-eigenvalues} {\rm If $\lambda_0$ is a real (multiple
or semi-simple) embedded eigenvalue of $T$, the Jordan chain can be
truncated at $f_n$ even if $[f_1,f_n] = 0$\footnote{The Fredholm
theory gives a necessary but not a sufficient condition for
existence of the solution $f_{n+1}$ in the Jordan chain
(\ref{Jordan-chain-formula}) if the eigenvalue $\lambda_0$ is
embedded into the continuous spectrum.}. In the latter case, the
neutral subspaces ${\cal H}_0$ for $n = 2k$ and $\tilde{\cal H}_0$
for $n = 2k+1$ in Lemma \ref{Jordan chain} do not have to be maximal
non-positive or non-negative subspaces, while the subspace ${\cal
H}_{\lambda_0}$ in Corollary \ref{positive(negative) combination}
could be sign-indefinite. The construction of a maximal non-positive
subspace for (multiple or semi-simple) embedded eigenvalues depend
on the computations of the projection matrix $[f_i,f_j]$ in the
eigenspace of eigenvectors ${\cal H}_{\lambda} = {\rm
span}\{f_1,...,f_n\}$. We shall simplify this unnecessary
complication by an assumption that {\em all embedded eigenvalues are
simple}, such that the corresponding one-dimensional eigenspace
${\cal H}_{\lambda_0}$ is either positive or negative or neutral
with respect to $[\cdot,\cdot]$. }
\end{Remark}

\begin{Remark}
\label{remark-sum-dimension} {\rm By Lemma \ref{Jordan chain} and 
Corollary \ref{positive(negative) combination}, if $\lambda_0$ is a 
real (multiple or semi-simple) isolated eigenvalue, then the sum of 
dimensions of the maximal non-positive and non-negative subspaces of 
${\cal H}_{\lambda_0}$ equals the dimension of ${\cal 
H}_{\lambda_0}$ (although the intersection of the two subspaces can 
be non-empty). By Remark \ref{embedded-eigenvalues}, the dimension 
of ${\cal H}_{\lambda_0}$ for a real embedded eigenvalue can however 
be smaller than the sum of dimensions of the maximal non-positive 
and non-negative subspaces. Therefore, the presence of embedded 
eigenvalues introduces a complication in the count of eigenvalues of 
the generalized eigenvalue problem 
(\ref{generalized_eigenvalue_problem}). This complication was 
neglected in the implicit count of embedded eigenvalues in 
\cite{KKS}.}
\end{Remark}

\begin{Lemma}
\label{zero quadratic form for compex} Let $\lambda_0 \in 
\mathbb{C}$, ${\rm Im}(\lambda_0) > 0$ be an eigenvalue of $T$ in 
$\Pi_{\kappa}$, ${\cal H}_{\lambda_0}$ be the corresponding 
eigenspace, and $\tilde{\cal H}_{\lambda_0} = \{ {\cal 
H}_{\lambda_0}, {\cal H}_{\bar{\lambda}_0}\} \subset \Pi_{\kappa}$. 
Then, the neutral subspace ${\cal H}_{\lambda_0}$ is the maximal 
sign-definite subspace of $\tilde{\cal H}_{\lambda_0}$, such that 
$[f,f] = 0$ $\forall f \in \mathcal H_{\lambda_0}$.
\end{Lemma}

\begin{Proof}
By Lemma \ref{ortspaces} with $\lambda = \mu = \lambda_0$, the 
eigenspace ${\cal H}_{\lambda_0}$ is orthogonal to itself with 
respect to $[\cdot,\cdot]$, such that ${\cal H}_{\lambda_0}$ is a 
neutral subspace of $\tilde{\cal H}_{\lambda_0}$. It remains to 
prove that ${\cal H}_{\lambda_0}$ is the maximal sign-definite 
subspace in $\tilde{\cal H}_{\lambda}$. Let ${\cal H}_{\lambda_0} =  
= {\rm span}\{f_1,f_2,...,f_n \}$, where $\{f_1,f_2,...,f_n\}$ is 
the Jordan chain of eigenvectors (\ref{Jordan-chain-formula}). 
Consider a subspace $\tilde{H}_0 =  {\rm span}\{f_1,f_2,...,f_n, 
\bar{f}_j \}$, where $1 \leq j \leq n$ and construct a linear 
combination of $f_{n+1-j}$ and $\bar{f}_j$:
$$
[f_{n+1-j} + \alpha \bar{f}_j,f_{n+1-j} + \alpha \bar{f}_j] = 2 {\rm 
Re}\left(\alpha [\bar{f}_j,f_{n+1-j}]\right), \qquad \alpha \in 
\mathbb{C}.
$$
By the Fredholm theory, $[f_n,\bar{f}_1] \neq 0$ and by virtue of 
(\ref{help-help}), $[\bar{f}_j,f_{n+1-j}] \neq 0$. As a result, the 
linear combination $f_{n+1-j} + \alpha \bar{f}_j$ is sign-indefinite 
with respect to $[\cdot,\cdot]$. 
\end{Proof}

\begin{Lemma}
\label{lemma-shift-zero} Let ${\cal H}_0$ be eigenspace of a
multiple zero eigenvalue of $K^{-1} A$ in ${\cal H}$ and
$\{f_1,...,f_n\}$ be the Jordan chain of eigenvectors, such that
$f_1 \in {\rm Ker}(A)$. Let $0 < \delta < |\sigma_{-1}|$, where
$\sigma_{-1}$ is the smallest negative eigenvalue of $K^{-1} A$. If
$\omega_+ > 0$ then $[f_1,f_n] = (Kf_1,f_n) \neq 0$, and
\begin{itemize}
\item If $n$ is odd, the subspace ${\cal H}_0$ corresponds to a positive eigenvalue of the
operator $(A + \delta K)$ if $[f_1,f_n] > 0$ and to a negative
eigenvalue if $[f_1,f_n] < 0$.

\item If $n$ is even, the subspace ${\cal H}_0$ corresponds to a positive eigenvalue of the
operator $(A + \delta K)$ if $[f_1,f_n] < 0$ and to a negative
eigenvalue if $[f_1,f_n] > 0$.
\end{itemize}
\end{Lemma}

\begin{Proof}
Let $\mu(\delta)$ be an eigenvalue of the operator $A + \delta K$
related to the subspace ${\cal H}_0$. By standard perturbation
theory for isolated eigenvalues \cite{Kato}, each eigenvalue
$\mu_j(\delta)$ is a continuous function of $\delta$ and
\begin{equation}
\label{limiting-relation} \lim_{\delta \to 0^+}
\frac{\mu(\delta)}{\delta^n} = (-1)^{n+1} \frac{(K f_1,
f_n)}{(f_1,f_1)}.
\end{equation}
If $\omega_+ > 0$, the zero eigenvalue of $A$ is isolated from the
essential spectrum of $K^{-1} A$, such that $[f_1,f_n] = (K f_1,f_n)
\neq 0$ by the Fredholm theory. The statement of the lemma follows 
from the limiting relation (\ref{limiting-relation}). Since no 
eigenvalues of $K^{-1} A$ exists in $(-\sigma_{-1},0)$, the 
eigenvalue $\mu(\delta)$ remains sign-definite for $0 < \delta < 
|\sigma_{-1}|$.
\end{Proof}

\begin{Remark}
\label{remark-shift-zero} {\rm The statement holds for the case
$\omega_+ = 0$ provided that $n = 1$, $[f_1,f_1] \neq 0$, and
$\mu(\delta) < \omega_{A + \delta K}$, where $\omega_{A + \delta K} 
> 0$ is defined below the decomposition 
(\ref{decomposition-A-K-delta}). }
\end{Remark}

\begin{Remark}
\label{remark-count-eigenvalues} {\rm We are concerned here in the
bounds on the numbers of eigenvalues in the generalized eigenvalue
problem (\ref{generalized_eigenvalue_problem}) in terms of the
numbers of isolated eigenvalues of self-adjoint operators $A$ and
$K^{-1}$. Let $N_p^-$ ($N_n^-$) be the number of negative
eigenvalues of the bounded or unbounded operator $K^{-1} A$ with the
account of their multiplicities whose eigenvectors are associated to
the maximal non-negative (non-positive) subspace of $\Pi_{\kappa}$ 
with respect to $[\cdot,\cdot] = (K \cdot,\cdot)$. Similarly, let 
$N_p^0$ ($N_n^0$) be the number of zero eigenvalues of $K^{-1} A$ 
with the account of their multiplicities and $N_p^+$ ($N_n^+$) be 
the number of positive eigenvalues of $K^{-1} A$ with the account of 
their multiplicities, such that the corresponding eigenvectors are 
associated to the maximal non-negative (non-positive) subspace of 
$\Pi_{\kappa}$. In the case of real isolated eigenvalues, the sum of 
dimensions of the maximal non-positive and non-negative subspaces 
equals the dimension of the subspace ${\cal H}_{\lambda_0}$ by 
Remark \ref{remark-sum-dimension}. The splitting of the dimension of 
${\cal H}_{\lambda_0}$ between $N_p(\lambda_0)$ and $N_n(\lambda_0)$ 
is obvious for each semi-simple isolated eigenvalue $\lambda_0$ in 
Corollary \ref{positive(negative) combination}. In the case of a 
multiple real isolated eigenvalue $\lambda_0$ of algebraic 
multiplicity $n$, the same splitting is prescribed in Lemma 
\ref{Jordan chain}:
\begin{enumerate}
\item If $n = 2k$, $k \in \mathbb{N}$, then $N_p(\lambda_0) = N_n(\lambda_0) = k$.
\item If $n = 2k+1$, $k \in \mathbb{N}$ and $[f_1,f_n] > 0$, then $N_p(\lambda_0) = k+1$ and $N_n(\lambda_0) = k$.
\item If $n = 2k+1$, $k \in \mathbb{N}$ and $[f_1,f_n] < 0$, then $N_p(\lambda_0) = k$ and $N_n(\lambda_0) = k+1$.
\end{enumerate}
In the case of a simple real embedded eigenvalue $\lambda_0$, the
numbers $N_p(\lambda_0)$ and $N_n(\lambda_0)$ for a one-dimensional
subspace ${\cal H}_{\lambda_0}$ are prescribed in Remark
\ref{embedded-eigenvalues} as follows:
\begin{enumerate}
\item If $[f_1,f_1] > 0$, then $N_p(\lambda_0) = 1$,
$N_n(\lambda_0) = 0$. \item If $[f_1,f_1] < 0$, then
$N_p(\lambda_0) = 0$, $N_n(\lambda_0) = 1$. \item If $[f_1,f_1] =
0$, then $N_p(\lambda_0) = N_n(\lambda_0) = 1$.
\end{enumerate}
In the case (iii), the sum $N_p(\lambda_0) + N_n(\lambda_0)$ exceeds
the dimension of ${\cal H}_{\lambda_0}$. Let $N_{c^+}$ ($N_{c^-}$)
be the number of complex eigenvalues in the upper half plane $\gamma
\in \mathbb{C}$, ${\rm Im}(\gamma) > 0$ (${\rm Im}(\gamma) < 0$).
The maximal sign-definite subspace of $\Pi_{\kappa}$ associated to
complex eigenvalues is prescribed by Lemma \ref{zero quadratic form 
for compex}.}
\end{Remark}

\begin{Theorem}
\label{equality 1} Let assumptions P1--P2 be satisfied. Eigenvalues
of the generalized eigenvalue problem
(\ref{generalized_eigenvalue_problem}) satisfy the pair of
equalities:
\begin{eqnarray}
\label{negative-index-A} N_p^- + N_n^0 + N_n^+ + N_{c^+} & = &  {\rm
dim}({\cal H}_{A+\delta K}^-) \\ \label{negative-index-K} N_n^- +
N_n^0 + N_n^+ + N_{c^+} & = & {\rm dim}({\cal H}_K^-)
\end{eqnarray}
\end{Theorem}

\begin{Proof}
We use the shifted eigenvalue problem
(\ref{shifted_eigenvalue_problem}) with sufficiently small $\delta
> 0$ and consider the bounded operator $T = (A + \delta K)^{-1}K$.
By Lemma \ref{BG product}, the operator $T$ is self-adjoint with
respect to $[\cdot,\cdot] = (K \cdot,\cdot)$ and it has a
non-positive invariant subspace of dimension $\kappa = {\rm
dim}({\cal H}_K^-)$. Counting all eigenvalues of the shifted
eigenvalue problem (\ref{shifted_eigenvalue_problem}) in Remark
\ref{remark-count-eigenvalues}, we establish the equality
(\ref{negative-index-K}). The other equality
(\ref{negative-index-A}) follows from a count for the bounded
operator $\tilde{T} = K (A + \delta K)^{-1}$ which is self-adjoint
with respect to $[\cdot,\cdot] = ((A + \delta K)^{-1}
\cdot,\cdot)$. The self-adjoint operator $(A + \delta K)^{-1}$
defines the indefinite metric in the Pontryagin space
$\tilde{\Pi}_{\tilde{\kappa}}$, where $\tilde{\kappa} = {\rm
dim}({\cal H}_{A + \delta K}^-)$. For any semi-simple eigenvalue
$\gamma_0$ of the shifted eigenvalue problem
(\ref{shifted_eigenvalue_problem}), we have
$$
\forall f,g \in {\cal H}_{\gamma_0}, \qquad ((A + \delta K) f,g) =
(\gamma_0 + \delta)(Kf,g).
$$
If $\gamma_0 \geq 0$ or ${\rm Im}(\gamma_0) \neq 0$, the maximal
non-positive eigenspace of $\tilde{T}$ in
$\tilde{\Pi}_{\tilde{\kappa}}$ associated with $\gamma_0$
coincides with the maximal non-positive eigenspace of $T$ in
$\Pi_{\kappa}$. If $\gamma_0 < 0$, the maximal non-positive
eigenspace of $\tilde{T}$ in $\tilde{\Pi}_{\tilde{\kappa}}$
coincides with the maximal non-negative eigenspace of $T$ in
$\Pi_{\kappa}$. The same statement can be proved for the case of
multiple isolated eigenvalues $\gamma_0$ when the eigenspace is
defined by the Jordan block of eigenvectors,
$$
A f_j = \gamma_0 K f_j + f_{j-1}, \qquad j = 1,...,n,
$$
where $f_0 = 0$. The dimension of the maximal non-positive
eigenspace of $\tilde{T}$ in $\tilde{\Pi}_{\tilde{\kappa}}$ is
then $N_p^- + N_n^0 + N_n^+ + N_{c^+}$.
\end{Proof}

\begin{Corollary}
Let $N_{\rm neg} = {\rm dim}({\cal H}_{A + \delta K}^-) + {\rm
dim}({\cal H}_K^-)$ be the total negative index of the shifted
generalized eigenvalue problem (\ref{shifted_eigenvalue_problem}).
Let $N_{\rm unst} = N_p^- + N_n^- + 2N_{c^+}$ be the total number of
unstable eigenvalues that include $N^- = N_p^- + N_n^-$ negative
eigenvalues $\gamma < 0$ and $N_c = N_{c^+} + N_{c^-}$ complex
eigenvalues with ${\rm Im}(\gamma) \neq 0$. Under the assumptions of
Theorem \ref{equality 1}, it is true that
\begin{equation}
\label{closure-negative-index} \Delta N = N_{\rm neg} - N_{\rm unst}
= 2 N_n^+ + 2 N_n^0 \geq 0.
\end{equation}
\end{Corollary}

\begin{Proof}
The equality (\ref{closure-negative-index}) follows by the sum of
(\ref{negative-index-A}) and (\ref{negative-index-K}).
\end{Proof}

\begin{Theorem}
\label{inequality positive} Let assumptions P1--P2 be satisfied and
$\omega_+ > 0$. Let $N_A = {\rm dim}({\cal H}_A^- \oplus {\cal
H}_A^0 \oplus {\cal H}_A^+)$ be the total number of isolated
eigenvalues of $A$. Let $N_K = {\rm dim}({\cal H}_K^- \oplus {\cal
H}_K^+)$ be the total number of isolated eigenvalues of $K$. Assume
that no embedded eigenvalues of the generalized eigenvalue problem
(\ref{generalized_eigenvalue_problem}) exist. Isolated eigenvalues
of the generalized eigenvalue problem
(\ref{generalized_eigenvalue_problem}) satisfy the inequality:
\begin{equation}
\label{upper-bound} N_p^- + N_p^0 + N_p^+ +  N_{c^+} \leq N_A + N_K.
\end{equation}
\end{Theorem}

\begin{Proof}
We prove this theorem by contradiction. Let $\Pi$ be a subspace in
$\Pi_{\kappa}$ spanned by eigenvectors of the generalized eigenvalue
problem (\ref{generalized_eigenvalue_problem}) which belong to
$N_p^-$ negative eigenvalues $\gamma < 0$, $N_p^0$ zero eigenvalues
$\gamma = 0$, $N_p^+$ positive isolated eigenvalues $0 < \gamma <
\omega_+ \omega_-$, and $N_{c^+}$ complex eigenvalues ${\rm
Im}(\gamma) > 0$. Let us assume that $N_p^- + N_p^0 + N_p^+ +
N_{c^+} > N_A + N_K$.  By Gram--Schmidt orthogonalization, there
exist a vector $h$ in $\Pi$ such that $(h,f) = 0$ and $(h,g) = 0$,
where $f \in {\cal H}_A^- \oplus {\cal H}_A^0 \oplus {\cal H}_A^+$
and $g \in {\cal H}_K^- \oplus {\cal H}_K^+$, such that $h \in
\mathcal{H}_{A}^{\sigma_e(A)} \cap \mathcal{H}_{K}^{\sigma_e(K)}$.
As a result,
$$
(Ah,h) \geq \omega_+ (h,h), \qquad (Kh,h) \leq \omega_-^{-1} (h,h),
$$
and
$$
(Ah,h) \geq \omega_+ \omega_- (Kh,h).
$$
On the other hand, since $h \in \Pi$, we represent $h$ by a linear
combination of the eigenvectors in the corresponding subspaces of
$\Pi$, such that
\begin{eqnarray*}
(A h, h) & = & \sum_{i,j} \alpha_i \bar{\alpha}_j (A h_i,h_j)\\
& = & \sum_{\gamma_i = \gamma_j < 0} \alpha_i \bar{\alpha}_j (A h_i,
h_j) + \sum_{\gamma_i = \gamma_j = 0} \alpha_i \bar{\alpha}_j (A
h_i, h_j) + \sum_{0 < \gamma_i = \gamma_j < \omega_+ \omega_-}
\alpha_i \bar{\alpha}_j  (A h_i, h_j),
\end{eqnarray*}
where we have used Lemma \ref{ortspaces} and Corollary \ref{zero
quadratic form for compex}. By Lemma \ref{Jordan chain}, the
non-zero values in $(Ah_i,h_j)$ occur only for eigenvalues with odd
algebraic multiplicity $n = 2k+1$ for eigenvectors $(A
f_{k+1},f_{k+1})$, where $f_{k+1}$ is the generalized eigenvector in
the basis for $\tilde{\cal H}_0$. Since all these cases are similar
to the case of simple eigenvalues, we obtain
\begin{eqnarray*}
(Ah,h) & = & \sum_{\gamma_j < 0} |\alpha_j|^2 (A h_j, h_j) +
\sum_{\gamma_j = 0} |\alpha_j|^2 (A h_j, h_j) + \sum_{0 < \gamma_j <
\omega_+ \omega_-} |\alpha_j|^2  (A h_j, h_j) \\ & = &
\sum_{\gamma_j < 0} \gamma_j |\alpha_j|^2 (K h_j, h_j) + \sum_{0 <
\gamma_j < \omega_+ \omega_-} \gamma_j |\alpha_j|^2  (K h_j, h_j) \\
& < & \omega_+ \omega_- \sum_{0 < \gamma_j < \omega_+ \omega_-}
|\alpha_j|^2 (K h_j, h_j),
\end{eqnarray*}
where we have used the fact that $(K h_j,h_j) \geq 0$ in $\Pi$. On
the other hand,
\begin{eqnarray*}
(K h, h) & = & \sum_{i,j} \alpha_i \bar{\alpha}_j (K h_i,h_j) \\
& = & \sum_{\gamma_j < 0} |\alpha_j|^2 (K h_j, h_j) + \sum_{\gamma_j
= 0} |\alpha_j|^2 (K h_j, h_j) + \sum_{0 < \gamma_j <
\omega_+ \omega_-} |\alpha_j|^2  (K h_j, h_j) \\
& \geq & \sum_{0 < \gamma_j < \omega_+ \omega_-} |\alpha_j|^2  (K
h_j, h_j).
\end{eqnarray*}
Therefore, $(Ah,h) < \omega_+ \omega_- (Kh,h)$, which is a
contradiction.
\end{Proof}

\begin{Corollary}
\label{final inequality} Let $N_{\rm total} = N_A + N_K$ be the
total number of isolated eigenvalues of operators $A$ and $K$. Let
$N_{\rm isol} = N_p^- + N_n^- + N^0_p + N^0_n + N^+_p + N^+_n +
N_{c^+} + N_{c^-}$ be the total number of isolated eigenvalues of
the generalized eigenvalue problem
(\ref{generalized_eigenvalue_problem}). Under the assumptions of
Theorem \ref{inequality positive}, it is true that
\begin{equation}
\label{total-number} N_{\rm isol} \leq N_{\rm total} + {\rm
dim}({\cal H}_K^-).
\end{equation}
\end{Corollary}

\begin{Proof}
The inequality (\ref{total-number}) follows by the sum of
(\ref{negative-index-K}) and (\ref{closure-negative-index}).
\end{Proof}

\begin{Remark}
{\rm If a simple embedded eigenvalue with the corresponding
eigenvector $f_1$ is included into consideration according to Remark
\ref{remark-count-eigenvalues}, then the left-hand-side of the
inequality (\ref{upper-bound}) is increased by the number of simple
embedded eigenvalues with $[f_1,f_1] \leq 0$. For each embedded
eigenvalue, the right-hand-side in the inequality
(\ref{total-number}) is reduced by two if $[f_1,f_1] \leq 0$ and it
remains the same if $[f_1,f_1] > 0$. }
\end{Remark}

\section{Applications}

We shall describe three applications of the general analysis which
are related to recent studies of stability of solitons and vortices
in the nonlinear Schr\"{o}dinger equation and solitons in the 
Korteweg--de Vries equation.

\subsection{Solitons of the scalar nonlinear Schr\"{o}dinger equation}

Consider a scalar nonlinear Schr\"{o}dinger (NLS) equation,
\begin{equation}
\label{NLS} i \psi_t = -\Delta \psi + F(|\psi|^2) \psi, \qquad
\Delta = \partial^2_{x_1 x_1} + ... + \partial^2_{x_d x_d},
\end{equation}
where $(x,t) \in \mathbb{R}^d \times \mathbb{R}$ and $\psi \in
\mathbb{C}$. For a suitable nonlinear function $F(|\psi|^2)$, where 
$F$ is $C^{\infty}$ and $F(0) = 0$, the NLS equation (\ref{NLS}) 
possesses a solitary wave solution $\psi = \phi(x) e^{i \omega t}$, 
where $\omega > 0$, $\phi : \mathbb{R}^d \to \mathbb{R}$, and 
$\phi(x)$ is exponentially decaying $C^\infty$ function. See 
\cite{M93} for existence and uniqueness of ground state solutions to 
the NLS equation (\ref{NLS}). Linearization of the NLS equation 
(\ref{NLS}) with the anzats,
\begin{equation}
\label{linearization} \psi = \left( \phi(x) + [u(x) + i w(x)]
e^{\lambda t} + [\bar{u}(x) + i \bar{w}(x)] e^{\bar{\lambda} t}
\right) e^{i \omega t},
\end{equation}
where $\lambda \in \mathbb{C}$ and $(u,w) \in \mathbb{C}^2$, results
in the linearized Hamiltonian problem (\ref{coupled-problem}), where
$L_{\pm}$ are Schr\"{o}dinger operators,
\begin{eqnarray}
\label{operator-L+} L_+ & = & - \Delta + \omega + F(\phi^2) + 2
\phi^2 F'(\phi^2), \\ \label{operator-L-} L_- & = & - \Delta +
\omega + F(\phi^2).
\end{eqnarray}
We note that $L_{\pm}$ are unbounded operators and
$\sigma_c(L_{\pm}) = [\omega,\infty)$ with $\omega_+ = \omega_- =
\omega > 0$. The kernel of $L_-$ includes at least one eigenvector
$\phi(x)$ and the kernel of $L_+$ includes at least $d$
eigenvectors $\partial_{x_j} \phi(x)$, $j = 1,...,d$. The Hilbert
space ${\cal X}$ is defined as ${\cal X} =
H^1(\mathbb{R}^d,\mathbb{C})$ and the main assumptions P1-P2 are
satisfied due to exponential decay of the potential functions
$F(\phi^2)$ and $\phi^2 F'(\phi^2)$. Theorems \ref{equality 1} and
\ref{inequality positive} give precise count of eigenvalues of the
stability problem $L_- L_+ u = - \lambda^2 u$, provided that the
numbers ${\rm dim}({\cal H}_K^-)$, ${\rm dim}({\cal H}_{A+\delta
K}^-)$, $N_K$ and $N_A$ can be computed from the count of isolated
eigenvalues of $L_{\pm}$. We shall illustrate these computations
with two examples.

\vspace{0.5cm}

{\bf Example 1.} Let $\phi(x)$ be the ground state solution such
that $\phi(x) > 0$ on $x \in \mathbb{R}^d$. By spectral theory,
${\rm Ker}(L_-) = \{ \phi(x) \}$ is one-dimensional and the subspace
${\cal H}_K^-$ in (\ref{decomposition-K}) is empty, such that the
corresponding Pontryagin space is $\Pi_0$ with $\kappa = 0$.

\begin{itemize}
\item By the equality (\ref{negative-index-K}), it follows
immediately that $N_n^- = N_n^0 = N_n^+ = N_{c^+} = 0$ such that the
spectrum of the problem (\ref{generalized_eigenvalue_problem}) is
real-valued and all eigenvalues are semi-simple.

\item Since eigenvectors of ${\rm Ker}(A)$ are in the positive
subspace of $K$, the number of zero eigenvalues of the problem
(\ref{generalized_eigenvalue_problem}) is given by $N_p^0 = {\rm
dim}({\cal H}^0_A)$. If $\frac{d}{d\omega} \| \phi \|^2_{L^2} =
0$, then $L_+
\partial_{\omega} \phi(x) = - \phi(x)$, $\partial_{\omega} \phi(x)
\in {\cal H}$, ${\cal P} \phi = 0$, and the eigenvector
$\partial_{\omega} \phi(x) \in {\rm Ker}({\cal P} L_+ {\cal P})$, 
such that $z_1 = 1$. If $z(L_+) = d$, then $z_0 = 0$ since 
$\partial_{x_j} \phi(x) \in {\cal H}$ for $j = 1,...,d$. In 
particular, $L_- x_j \phi(x) = - 2
\partial_{x_j} \phi(x)$ and $(K \partial_{x_j} \phi, \partial_{x_j}
\phi) = \frac{1}{4} \| \phi \|^2_{L^2} > 0$.

\item By the equality (\ref{negative-index-A}), the number of
negative eigenvalues of the problem
(\ref{generalized_eigenvalue_problem}) is given by $N_p^- = {\rm
dim}({\cal H}_{A + \delta K}^-)$. By Lemma \ref{lemma-shift-zero}
with $n = 1$ and $[f_1,f_1] > 0$, all zero eigenvalues of $A$ become 
positive eigenvalues of $A + \delta K$ for $\delta > 0$. By 
Propositions \ref{lemma-constrained-space} and 
\ref{proposition-zero-splitting}, we have $ {\rm dim}({\cal H}_{A + 
\delta K}^-) = {\rm dim}({\cal H}_A^-) = n(L_+) - n_0$. By Theorem 
3.1 in \cite{GSS2}, $n_0 = 1$ if $\frac{d}{d\omega} \| \phi 
\|^2_{L^2} > 0$ and $n_0 = 0$ otherwise.

\item By the inequality (\ref{upper-bound}), the number of
positive eigenvalues of the problem
(\ref{generalized_eigenvalue_problem}) are bounded from above by
$N_p^+ \leq {\rm dim}({\cal H}^+_A) + {\rm dim}({\cal H}^+_K)$. By
Proposition \ref{lemma-constrained-space}, it is then $N_p^+ \leq
p(L_+) + p(L_-) + n_0 + z_0 - z_1$.
\end{itemize}

\begin{Remark}
{\rm Under the assumptions that $n(L_+) = 1$, $z(L_+) = d$,
$p(L_+) = p(L_-) = 0$ and $\frac{d}{d\omega} \| \phi \|^2_{L^2} <
0$, it follows from the above properties that $N_p^- = 1$, $N_p^0
= d$, and $N_p^+ = 0$.  The assumptions $n(L_+) = 1$, $z(L_+) = d$
and $\frac{d}{d\omega} \| \phi \|^2_{L^2} < 0$ are verified for
the super-critical NLS equation with the power nonlinearity $F =
|\psi|^p$ \cite{Schlag}. Under the further assumptions $p(L_+) =
p(L_-) = 0$, the statement $N_p^+ = 0$ is proved directly in
Proposition 2.1.2 \cite{Per} and Proposition 9.2 \cite{KSchlag}
for $d = 1$ and in Lemma 1.8 \cite{Schlag} for $d = 3$. The same
statement follows here by the count of eigenvalues from Theorem
\ref{inequality positive}. }
\end{Remark}

\begin{Remark}
{\rm Systems of coupled NLS equations generalize the scalar NLS
equation (\ref{NLS}). Stability of vector solitons in the coupled
NLS equations results in the same linearized Hamiltonian system
(\ref{coupled-problem}) with the matrix Schr\"{o}dinger operators
$L_{\pm}$. General results for non-ground state solutions are
obtained in \cite{KKS,Pel} for $d = 1$ and in \cite{CPV} for $d =
3$. Multiple and embedded eigenvalues were either excluded from
analysis by an assumption \cite{Pel,CPV} or were treated implicitly
\cite{KKS}. The present manuscript generalizes these results for an
abstract case with a precise count of multiple and embedded
eigenvalues.}
\end{Remark}

{\bf Example 2.} Let the scalar NLS equation (\ref{NLS}) with $F =
|\psi|^2$ be discretized with $\Delta \equiv \epsilon \Delta_{\rm
disc}$, where $\Delta_{\rm disc}$ is the second-order discrete
Laplacian and $\epsilon$ is a small parameter. We note that
$\Delta_{\rm disc}$ is a bounded operator and
$\sigma_c(-\Delta_{\rm disc}) \in [0,4d]$. By the
Lyapunov--Schmidt reduction method, the solution $\psi = \phi e^{i
\omega t}$ with $\omega
> 0$ bifurcates from the limit $\epsilon = 0$ with $N$ non-zero
lattice nodes to $\epsilon \neq 0$, such that $\frac{d}{d\omega}
\| \phi \|^2_{l^2} > 0$, ${\rm ker}(L_+) = \varnothing$, and ${\rm
ker}(L_-) = \{ \phi \}$ is one-dimensional for $0 < |\epsilon| <
\epsilon_0$ with $\epsilon_0
> 0$ (see \cite{PKF1} for $d = 1$ and \cite{PKF2} for $d = 2$). By
the equalities (\ref{negative-index-A}) and
(\ref{negative-index-K}), we have
\begin{eqnarray*}
N_p^- + N_n^+ + N_{c^+} & = &  n(L_+) - 1, \\ N_n^- + N_n^+ +
N_{c^+} & = & n(L_-),
\end{eqnarray*}
where $n(L_+) = N$ and $n(L_-) \leq N - 1$ in the domain $0 < 
|\epsilon | < \epsilon_0$. When small positive eigenvalues of $L_-$ 
are simple for $\epsilon \neq 0$, it is true that $N_n^+ = n(L_-)$, 
$N_n^- = N_{c^+} = 0$, and $N_p^- = N - 1 - n(L_-)$ (see Corollary 
3.5 in \cite{PKF1}\footnote{Corollary 3.5 in \cite{PKF1} is valid 
only when small positive eigenvalues of $L_-$ are simple. It is 
shown in \cite{PKF2} that the case of multiple small positive 
eigenvalues of $L_-$ leads to splitting of real eigenvalues $N_p^-$ 
of the generalized eigenvalue problem 
(\ref{generalized_eigenvalue_problem}) to complex eigenvalues 
$N_{c^+}$ beyond the leading-order Lyapunov--Schmidt reduction. The 
case of multiple small negative eigenvalues of $L_-$ does not lead 
to this complication since the semi-simple purely imaginary 
eigenvalues $N_n^+$ do not split to complex eigenvalues 
$N_{c^+}$.}). It is clear that the Lyapunov--Schmidt reduction 
method gives a more precise information on numbers $N_n^{\pm}$, 
$N_p^-$, and $N_{c^+}$, compared to the general equalities above. 
Similarly, by the inequality (\ref{upper-bound}) and the counts 
above for $0 < |\epsilon| < \epsilon_0$, we have
$$
N_p^+ \leq 2 n(L_-) + {\rm dim}({\cal H}_A^+) + {\rm dim}({\cal
H}_K^+).
$$
When $0 < |\epsilon| < \epsilon_0$, the numbers $N_p^+$ and ${\rm 
dim}({\cal H}_A^+)$ give precisely the numbers of edge bifurcations 
from the essential spectrum of $K^{-1} A$ and $A$ respectively, 
while the number ${\rm dim}({\cal H}_K^+)$ exceeds the number of 
edge bifurcations from the essential spectrum of $K^{-1}$ by $N - 1 
- n(L_-)$. When the solution $\phi$ is a ground state, we have $N = 
1$ and $n(L_-) = 0$, such that the edge bifurcations in the spectrum 
of $K^{-1}A$ may only occur if there are edge bifurcations in the 
spectrum of $A$ and $K^{-1}$ for $\epsilon \neq 0$. When $N > 1$, 
the relation between edge bifurcations in the self-adjoint and 
non-self-adjoint problems become less direct.

\begin{Remark}
{\rm The Lyapunov--Schmidt reduction method was also used for
continuous coupled NLS equations with and without external
potentials. See \cite{kap1,kap2,pelyang} for various results on
the count of unstable eigenvalues in parameter continuations of
the NLS equations. }
\end{Remark}

\subsection{Vortices of the scalar nonlinear Schr\"{o}dinger equation}

Consider the scalar two-dimensional NLS equation (\ref{NLS}) with
$d = 2$ in polar coordinates $(r,\theta)$:
\begin{equation}
\label{NLSrad} i \psi_t = -\Delta \psi + F(|\psi|^2) \psi, \qquad
\Delta = \partial^2_{r r} + \frac{1}{r} \partial_r + \frac{1}{r^2}
\partial^2_{\theta \theta},
\end{equation}
where $r > 0$ and $\theta \in [0,2\pi]$. Assume that the NLS
equation (\ref{NLSrad}) possesses a charge-$m$ vortex solution $\psi 
= \phi(r) e^{i m \theta + i \omega t}$, where $\omega > 0$, $m \in 
\mathbb{N}$, $\phi  : \mathbb{R}_+ \to \mathbb{R}$, and $\phi(r)$ is 
exponentially decaying $C^\infty$ function for $r > 0$ with $\phi(0) 
= 0$. See \cite{PegWar} for existence results of charge-$m$ vortices 
in the cubic-quintic NLS equation with $F = - |\psi|^2 + |\psi|^4$. 
Linearization of the NLS equation (\ref{NLSrad}) with the anzats,
\begin{equation}
\label{linearization_rad} \psi = \left( \phi(r) e^{i m \theta} +
\varphi_+(r,\theta) e^{\lambda t} + \bar{\varphi}_-(r,\theta)
e^{\bar{\lambda} t} \right) e^{i \omega t},
\end{equation}
where $\lambda \in \mathbb{C}$ and $(\varphi_+,\varphi_-) \in
\mathbb{C}^2$, results in the stability problem,
\begin{equation}
\label{vortex-stability} \sigma_3 H \mbox{\boldmath $\varphi$} = i
\lambda \mbox{\boldmath $\varphi$},
\end{equation}
where $\mbox{\boldmath $\varphi$} = (\varphi_+,\varphi_-)^T$,
$\sigma_3 = {\rm diag}(1,-1)$, and
$$
H = \left( \begin{array}{cc} - \Delta + \omega + F(\phi^2) + \phi^2
F'(\phi^2) & \phi^2 F'(\phi^2) e^{2 i m \theta} \\ \phi^2 F'(\phi^2)
e^{-2 i m \theta} & - \Delta + \omega + F(\phi^2) + \phi^2
F'(\phi^2) \end{array} \right).
$$
Expand $\mbox{\boldmath $\varphi$}(r,\theta)$ in the Fourier series
$$
\mbox{\boldmath $\varphi$} = \sum_{n \in \mathbb{Z}} \mbox{\boldmath
$\varphi$}^{(n)}(r) e^{ i n \theta}
$$
and reduce the problem to a sequence of spectral problems for ODEs:
\begin{equation}
\label{stability-vortex} \sigma_3 H_n \mbox{\boldmath $\varphi$}_n =
i \lambda \mbox{\boldmath $\varphi$}_n, \qquad n \in \mathbb{Z},
\end{equation}
where $\mbox{\boldmath $\varphi$}_n =
(\varphi_+^{(n+m)},\varphi_-^{(n-m)})^T$, and
{\small
$$ H_n =
\left( \begin{array}{cc} - \partial^2_{rr} - \frac{1}{r}
\partial_r + \frac{(n+m)^2}{r^2} + \omega + F(\phi^2) + \phi^2
F'(\phi^2) & \phi^2 F'(\phi^2) \\ \phi^2 F'(\phi^2) & -
\partial^2_{rr} - \frac{1}{r} \partial_r + \frac{(n-m)^2}{r^2} + \omega
+ F(\phi^2) + \phi^2 F'(\phi^2) \end{array} \right).
$$
} When $n = 0$, the stability problem (\ref{stability-vortex})
transforms to the linearized Hamiltonian system
(\ref{coupled-problem}), where $L_{\pm}$ is given by
(\ref{operator-L+})--(\ref{operator-L-}) with $\Delta =
\partial^2_{r r} + \frac{1}{r} \partial_r - \frac{m^2}{r^2}$ and $(u,w)$ are
given by $u = \varphi_+^{(m)} + \varphi_-^{(-m)}$ and $w =
-i(\varphi_+^{(m)} - \varphi_-^{(-m)})$. When $n \in \mathbb{N}$,
the stability problem (\ref{stability-vortex}) transforms to the
linearized Hamiltonian system (\ref{coupled-problem}) with $L_+ =
H_n$ and $L_- = \sigma_3 H_n \sigma_3$, where
$$
L_+ = L_- + 2 \phi^2 F'(\phi^2) \sigma_1, \qquad \sigma_1 = \left(
\begin{array}{cc} 0 & 1 \\ 1 & 0 \end{array} \right),
$$
and $(u,w)$ are given by $u = \mbox{\boldmath $\varphi$}_n$ and $w
= - i \sigma_3 \mbox{\boldmath $\varphi$}_n$. When $-n \in
\mathbb{N}$, the spectrum of the stability problem
(\ref{stability-vortex}) can be obtained from that for $n \in
\mathbb{N}$ by the correspondence $H_{-n} = \sigma_1 H_n
\sigma_1$.

Let us introduce the weighted inner product for functions on $r \geq 
0$:
$$
(f,g) = \int_0^{\infty} f(r) g(r) r dr.
$$
In all cases $n = 0$, $n \in \mathbb{N}$ and $-n \in \mathbb{N}$,
$L_{\pm}$ are unbounded self-adjoint differential operators and
$\sigma_c(L_{\pm}) = [\omega,\infty)$, such that $\omega_+ =
\omega_- = \omega > 0$. The kernel of $H_n$ includes at least one
eigenvector for $n = \pm 1$
$$
\mbox{\boldmath $\phi$}_{\pm 1} = \phi'(r) {\bf 1} \mp \frac{m}{r} 
\phi(r) \sigma_3 {\bf 1}, \qquad {\bf 1} = (1,1)^T,
$$
and at least one eigenvector for $n = 0$: $\mbox{\boldmath $\phi$}_0 
= \phi(r) \sigma_3 {\bf 1}$. The Hilbert space ${\cal X}$ associated 
with the weighted inner product is defined as ${\cal X} = 
H^1(\mathbb{R}_+,\mathbb{C})$ for $n = 0$ and ${\cal X} = 
H^1(\mathbb{R}_+,\mathbb{C}^2)$ for $\pm n \in \mathbb{N}$. In all 
cases, the main assumptions P1-P2 are satisfied due to exponential 
decay of the potential functions $F(\phi^2)$ and $\phi^2 
F'(\phi^2)$.

The case $n = 0$ is the same as for solitons (see Section 5.1). We
shall hence consider adjustments in the count of eigenvalues in
the case $\pm n \in \mathbb{N}$, when the stability problem
(\ref{stability-vortex}) is rewritten in the form,
\begin{equation}
\label{stability-vortex-2} \left\{ \begin{array}{c} \sigma_3 H_n
\mbox{\boldmath $\varphi$}_n = i \lambda \mbox{\boldmath
$\varphi$}_n \\
\sigma_3 H_{-n} \mbox{\boldmath $\varphi$}_{-n} = i \lambda
\mbox{\boldmath $\varphi$}_{-n}\end{array} \right. \qquad n \in
\mathbb{N}.
\end{equation}
Let $L_+$ be a diagonal composition of $H_n$ and $H_{-n}$ and
$L_-$ be a diagonal composition of $\sigma_3 H_n \sigma_3$ and
$\sigma_3 H_{-n} \sigma_3$.

\begin{Lemma}
\label{lemma-stability-vortex} Let $\lambda$ be an eigenvalue of
the stability problem (\ref{stability-vortex-2}) with the
eigenvector $(\mbox{\boldmath $\varphi$}_n,{\bf 0})$. Then there
exists another eigenvalue $-\lambda$ with the linearly independent
eigenvector $({\bf 0},\sigma_1 \mbox{\boldmath $\varphi$}_n)$. If
${\rm Re}(\lambda)
> 0$, there exist two more eigenvalues
$\bar{\lambda}$,$-\bar{\lambda}$ with the linearly independent
eigenvectors  $({\bf 0},\sigma_1 \bar{\mbox{\boldmath
$\varphi$}}_n)$, $(\bar{\mbox{\boldmath $\varphi$}}_n,{\bf 0})$.
\end{Lemma}

\begin{Proof}
We note that $\sigma_1 \sigma_3 = - \sigma_3 \sigma_1$ and
$\sigma_1^2 = \sigma_3^2 = \sigma_0$, where $\sigma_0 = {\rm
diag}(1,1)$. Therefore, each eigenvalue $\lambda$ of $H_n$ with
the eigenvector $\mbox{\boldmath $\varphi$}_n$ generates
eigenvalue $-\lambda$ of $H_{-n}$ with the eigenvector
$\mbox{\boldmath $\varphi$}_{-n} = \sigma_1 \mbox{\boldmath
$\varphi$}_n$. When ${\rm Re}(\lambda) \neq 0$, each eigenvalue
$\lambda$ of $H_n$ generates also eigenvalue $-\bar{\lambda}$ of
$H_n$ with the eigenvector $\bar{\mbox{\boldmath $\varphi$}}_n$
and eigenvalue $\bar{\lambda}$ of $H_{-n}$ with the eigenvector
$\mbox{\boldmath $\varphi$}_{-n} = \sigma_1 \bar{\mbox{\boldmath
$\varphi$}}_n$.
\end{Proof}

\begin{Theorem}
\label{theorem-vortex} Let $N_{\rm real}$ be the number of real
eigenvalues in the stability problem (\ref{stability-vortex-2})
with ${\rm Re}(\lambda) > 0$, $N_{\rm comp}$ be the number of
complex eigenvalues with ${\rm Re}(\lambda) > 0$ and ${\rm
Im}(\lambda) > 0$, $N_{\rm imag}^-$ be the number of purely
imaginary eigenvalues with ${\rm Im}(\lambda) > 0$ and
$(\mbox{\boldmath $\varphi$}_n,H_n \mbox{\boldmath $\varphi$}_n)
\leq 0$, and $N_{\rm zero}^-$ is the algebraic multiplicity of the
zero eigenvalue of $\sigma_3 H_n \mbox{\boldmath $\varphi$}_n = i
\lambda \mbox{\boldmath $\varphi$}_n$ counted in Remark
\ref{remark-count-eigenvalues}. Assume that all embedded
eigenvalues are simple. Then,
\begin{equation}
\label{closure-relation2} \frac{1}{2} N_{\rm real} + N_{\rm comp}
= n(H_n) - N_{\rm zero}^- - N_{\rm imag}^-,
\end{equation}
where $N_{\rm real}$ is even. Moreover, $n_0 = n_-$, where $n_0$
and $n_-$ are defined by Propositions
\ref{lemma-constrained-space} and
\ref{proposition-zero-splitting}.
\end{Theorem}

\begin{Proof}
By Lemma \ref{lemma-stability-vortex}, a pair of real eigenvalues
of $\sigma_3 H_n \mbox{\boldmath $\varphi$}_n = i \lambda
\mbox{\boldmath $\varphi$}_n$ corresponds to two linearly
independent eigenvectors $\mbox{\boldmath $\varphi$}_n$ and
$\bar{\mbox{\boldmath $\varphi$}}_n$. Because $(H_n
\mbox{\boldmath $\varphi$}_n,\mbox{\boldmath $\varphi$}_n)$ is
real-valued and hence zero for $\lambda \in \mathbb{R}$, we have
$$
\left( H_n (\mbox{\boldmath $\varphi$}_n \pm \bar{\mbox{\boldmath
$\varphi$}}_n),(\mbox{\boldmath $\varphi$}_n \pm
\bar{\mbox{\boldmath $\varphi$}}_n)\right) = \pm 2 {\rm Re}(H_n
\mbox{\boldmath $\varphi$}_n, \bar{\mbox{\boldmath $\varphi$}}_n).
$$
By counting multiplicities of the real negative and complex
eigenvalues of the problem (\ref{generalized_eigenvalue_problem})
associated to the stability problem (\ref{stability-vortex-2}), we
have $N_n^- = N_p^- = N_{\rm real}$ and $N_{c^+} = 2 N_{\rm
comp}$. By Lemma \ref{lemma-stability-vortex}, a pair of purely
imaginary and zero eigenvalues of the stability problem
(\ref{stability-vortex-2}) corresponds to two linearly independent
eigenvectors $(\mbox{\boldmath $\varphi$}_n,{\bf 0})$ and $({\bf
0}, \mbox{\boldmath $\varphi$}_{-n})$, where $ \mbox{\boldmath
$\varphi$}_{-n} = \sigma_1 \mbox{\boldmath $\varphi$}_n$ and
$(H_{-n} \mbox{\boldmath $\varphi$}_{-n}, \mbox{\boldmath
$\varphi$}_{-n}) = (H_n \mbox{\boldmath $\varphi$}_n,
\mbox{\boldmath $\varphi$}_n)$. By counting multiplicities of the
real positive and zero eigenvalues of the problem
(\ref{generalized_eigenvalue_problem}) associated to the stability
problem (\ref{stability-vortex-2}), we have $N_n^0 = 2 N_{\rm
zero}^-$ and $N_n^+ = 2 N_{\rm imag}^-$. Since the spectra of
$H_n$, $\sigma_1 H_n \sigma_1$, and $\sigma_3 H_n \sigma_3$
coincide, we have $n(L_+) = n(L_-) = 2n(H_n)$. As a result, the
equality (\ref{closure-relation2}) follows by the equality
(\ref{negative-index-K}) of Theorem \ref{equality 1}. By Lemma
\ref{lemma-stability-vortex}, the multiplicity of $N_{\rm real}$
is even in the stability problem (\ref{stability-vortex-2}). The
other equality (\ref{negative-index-A}) of Theorem \ref{equality
1} recovers the same answer provided that ${\rm dim}({\cal H}_{A +
\delta K}^-) = {\rm dim}({\cal H}_A^-) + n_- = n(L_+) - n_0 + n_-
= n(L_+)$, such that $n_0 = n_-$.
\end{Proof}

{\bf Example 3.} Let $\phi(r)$ be the fundamental charge-$m$
vortex solution such that $\phi(r) > 0$ for $r > 0$ and $\phi(0) =
0$. By spectral theory, ${\rm Ker}(H_0)$ is one-dimensional with
the eigenvector $\mbox{\boldmath $\phi$}_0$. The analysis of $n =
0$ is similar to Example 1. In the case $n \in \mathbb{N}$, we
shall assume that ${\rm Ker}(H_1) = \{ \mbox{\boldmath $\phi$}_1
\}$ and ${\rm Ker}(H_n) = \varnothing$ for $n \geq 2$.

\begin{itemize}
\item Since $(\sigma_3 \mbox{\boldmath $\phi$}_1,\mbox{\boldmath
$\phi$}_1) = 0$ and ${\rm Ker}(\sigma_3 H_1 \sigma_3) = \{
\sigma_3 \mbox{\boldmath $\phi$}_1 \}$, then $\mbox{\boldmath
$\phi$}_1 \in {\cal H}$, such that $z_0 = 0$.

\item By direct computation, $(\sigma_3 H_1 \sigma_3)^{-1}
\mbox{\boldmath $\phi$}_1 = -\frac{1}{2} r \phi(r) {\bf 1}$ and
$$
((\sigma_3 H_1 \sigma_3)^{-1} \mbox{\boldmath
$\phi$}_1,\mbox{\boldmath $\phi$}_1) = \int_0^{\infty} r \phi^2(r)
dr > 0.
$$
By Lemma \ref{lemma-shift-zero}, we have $N_n^0 = n_- = 0$ for $n
= 1$, as well as for $n \geq 2$. By Proposition
\ref{lemma-constrained-space}, we have $A(0) < 0$ such that $n_0 =
z_1 = 0$ for all $n \in \mathbb{N}$.

\item By Theorem \ref{theorem-vortex}, we have
\begin{equation}
\label{closure-10} N_{\rm real} + 2 N_{\rm comp} = 2 n(H_n) - 2
N_{\rm imag}^-.
\end{equation}
If all purely imaginary eigenvalues are semi-simple and isolated,
$N_{\rm imag}^-$ gives the total number of eigenvalues in the
stability problem (\ref{stability-vortex-2}) with ${\rm
Re}(\lambda) = 0$, ${\rm Im}(\lambda) > 0$, and negative Krein
signature $(H_n \mbox{\boldmath $\varphi$}_n,\mbox{\boldmath
$\varphi$}_n) < 0$.
\end{itemize}

\begin{Remark}
{\rm Stability of vortices was considered numerically in
\cite{PegWar}, where Lemma \ref{lemma-stability-vortex} was also
obtained. The closure relation (\ref{closure-10}) was also discussed 
in \cite{KKS} in a more general context. Detailed comparison of 
numerical results and the closure relation (\ref{closure-10}) can be 
found in \cite{Kollar}.  Vortices in the discretized scalar NLS 
equation were considered with the Lyapunov--Schmidt reduction method 
in \cite{PKF2}. Although the reduced eigenvalue problems were found 
in a much more complicated form compared to the reduced eigenvalue 
problem for solitons, the relation (\ref{closure-10}) was confirmed 
in all particular vortex configurations considered in \cite{PKF2}.  
}
\end{Remark}

\begin{Remark}
{\rm  Since $N_{\rm real}$ is even, splitting of simple complex 
eigenvalues $N_{\rm comp}$ into double real eigenvalues $N_{\rm 
real}$ is prohibited by the closure relation (\ref{closure-10}). 
Simple complex eigenvalues may either coalesce into a double real 
eigenvalue that persists or reappear again as simple complex 
eigenvalues. }
\end{Remark}

\subsection{Solitons of the fifth-order Korteweg--De Vries equation}

Consider a general fifth-order KdV equation,
\begin{equation}
\label{general-KdV} v_t = a_1 v_x - a_2 v_{xxx} + a_3 v_{xxxxx} +
3 b_1 v v_x - b_2 \left( v v_{xxx} + 2 v_x v_{xx} \right) + 6 b_3
v^2 v_x,
\end{equation}
where $(a_1,a_2,a_3)$ and $(b_1,b_2,b_3)$ are real-valued
coefficients for linear and nonlinear terms, respectively. Without
loss of generality, we assume that $a_3 > 0$ and
\begin{equation}
\label{wave-speed} c_{\rm wave}(k) = a_1 + a_2 k^2 + a_3 k^4 \geq
0, \qquad \forall k \in \mathbb{R}.
\end{equation}
For suitable values of parameters, there exists a traveling wave
solution $v(x,t)= \phi(x-ct)$, where $c > 0$ and $\phi : \mathbb{R}
\mapsto \mathbb{R}$, such that the function $\phi(x)$ is even and
exponentially decaying as $|x| \to \infty$. Existence of traveling
waves was established in \cite{Z87,HS88,AT92} for $b_2 = b_3 = 0$,
in \cite{CG97} for $b_3 = 0$, in \cite{K97} for $b_1 = -b_2 = b_3 =
1$, and in \cite{L99} for $b_3 = 0$ or $b_1 = b_2= 0$. Linearization
of the fifth-order KdV equation (\ref{general-KdV}) with the anzats
$$
v(x,t) = \phi(x-ct) + w(x-ct) e^{\lambda t}
$$
results in the stability problem
\begin{equation}
\label{spectrum}
\partial_x L_- w = \lambda w,
\end{equation}
where $L_-$ is an unbounded fourth-order operator,
\begin{equation}
\label{operator-L} L_- = a_3 \frac{d^4}{d x^4} - a_2 \frac{d^2}{d
x^2} + a_1 + c + 3 b_1 \phi(x) - b_2 \frac{d}{d x} \phi(x)
\frac{d}{d x} - b_2 \phi''(x) + 6 b_3 \phi^2(x).
\end{equation}
With the account of the condition (\ref{wave-speed}), the continuous 
spectrum of $L_-$ is located for $\sigma_c(L_-) \in [c,\infty)$, 
such that $\omega_- = c > 0$. The kernel of $L_-$ includes at least 
one eigenvector $\phi'(x)$. Since the image of $L_-$ is in 
$L^2(\mathbb{R})$, the eigenfunction $w(x) \in L^1(\mathbb{R})$ for 
$\lambda \neq 0$ satisfies the constraint:
\begin{equation}
\label{mass-constraint} (1, w) = \int_{\mathbb{R}} w(x) dx = 0.
\end{equation}
Let $w = u'(x)$, where $u(x) \to 0$ as $|x| \to \infty$ and define
$L_+ = - \partial_x L_- \partial_x$. The continuous spectrum of
$L_+$ is located for $\sigma_c(L_+) \in [0,\infty)$, such that
$\omega_+ = 0$. The kernel of $L_+$ includes at least one
eigenvector $\phi(x)$.

Let the Hilbert space ${\cal X}$ be defined as ${\cal X} =
H^3(\mathbb{R},\mathbb{C})$. The main assumptions P1-P2 for $L_-$
and $L_+$ are satisfied due to exponential decay of the potential
function $\phi(x)$. The stability problem (\ref{spectrum}) is
equivalent to the linearized Hamiltonian system
(\ref{coupled-problem}). The operator $L_-$ defines the constrained 
subspace ${\cal H}$ and hence the Pontryagin space $\Pi_{\kappa}$ 
with $\kappa = n(L_-)$. The operator $L_+$ has the embedded kernel 
to the endpoint of the essential spectrum of $L_+$. This introduces 
a technical complication in computations of the inverse of $L_+$ 
\cite{KP}, which we avoid here with the use of the shifted 
generalized eigenvalue problem (\ref{shifted_eigenvalue_problem}) 
with $\delta > 0$. It is easy to estimate that
$$
\omega_{A + \delta K} = \inf_{k \in \mathbb{R}} \left[ k^2( c +
c_{\rm wave}(k)) + \frac{\delta}{c + c_{\rm wave}(k)} \right] \geq
\frac{\delta}{c} > 0,
$$
where $\omega_{A + \delta K}$ is defined below the decomposition
(\ref{decomposition-A-K-delta}). Theorem \ref{equality 1} can be
applied after appropriate adjustments in the count of isolated and
embedded eigenvalues in the stability problem (\ref{spectrum}).
Since $\omega_+ = 0$, Theorem \ref{inequality positive} is not
applicable to the fifth-order KdV equation (\ref{general-KdV}).
Because the continuous spectrum of $\partial_x L_-$ is on $\lambda
\in i \mathbb{R}$, all real and complex eigenvalues are isolated and
all purely imaginary eigenvalues including the zero eigenvalue are
embedded.

\begin{Lemma}
\label{lemma-inversion} Let $\lambda_j$ be a real eigenvalue of the
stability problem (\ref{spectrum}) with the real-valued eigenvector
$w_j(x)$, such that ${\rm Re}(\lambda_j) > 0$ and ${\rm
Im}(\lambda_j) = 0$. Then there exists another eigenvalue
$-\lambda_j$ in the problem (\ref{spectrum}) with the linearly
independent eigenvector $w_j(-x)$. The linear combinations
$w_j^{\pm}(x) = w_j(x) \pm w_j(-x)$ are orthogonal with respect to
the operator $L_-$,
\begin{eqnarray}
\label{orthogonality2} \left( L_- w_j^{\pm}, w_j^{\pm} \right) =
\pm 2 \rho_j, \qquad \left( L_- w_j^{\mp}, w_j^{\pm} \right) = 0,
\end{eqnarray}
where $\rho_j = \left( L_- w_j(-x),w_j(x)\right)$.
\end{Lemma}

\begin{Proof}
Since $\phi(-x) = \phi(x)$, the self-adjoint operator $L_-$ is
invariant with respect to the transformation $x \mapsto -x$. The
functions $w_j(x)$ and $w_j(-x)$ are linearly independent since
$w_j(x)$ has both symmetric and anti-symmetric parts provided that
$\lambda_j \neq 0$. Under the same constraint,
$$
\left( L_- w_j(\pm x), w_j(\pm x) \right) = \pm \lambda_j^{-1}
\left( L_- w_j(\pm x), \partial_x L_- w_j(\pm x) \right) = 0,
$$
and the orthogonality relations (\ref{orthogonality2}) hold by
direct computations.
\end{Proof}

\begin{Corollary}
\label{remark-complex-KdV} Let $\lambda_j$ be a complex eigenvalue
of the stability problem (\ref{spectrum}) with the complex-valued
eigenvector $w_j(x)$, such that ${\rm Re}(\lambda_j)
> 0$ and ${\rm Im}(\lambda_j) > 0$. Then there exist eigenvalues
$\bar{\lambda}_j$, $-\lambda_j$, and $-\bar{\lambda}_j$ in the
problem (\ref{spectrum}) with the linearly independent
eigenvectors $\bar{w}_j(x)$, $w_j(-x)$, and $\bar{w}_j(-x)$,
respectively.
\end{Corollary}

\begin{Lemma}
\label{lemma-orthogonality3} Let $\lambda_j$ be an embedded
eigenvalue of the stability problem (\ref{spectrum}) with the
complex-valued eigenvector $w_j(x)$, such that ${\rm Re}(\lambda_j)
= 0$ and ${\rm Im}(\lambda_j) > 0$. Then there exists another
eigenvalue $-\lambda_j = \bar{\lambda}_j$ in the problem
(\ref{spectrum}) with the linearly independent eigenvector $w_j(-x)
= \bar{w}_j(x)$. The linear combinations $w_j^{\pm}(x) = w_j(x) \pm
\bar{w}_j(x)$ are orthogonal with respect to the operator $L_-$,
\begin{eqnarray}
\label{orthogonality5} \left( L_- w_j^{\pm}, w_j^{\pm} \right) = 2
\rho_j, \qquad \left( L_- w_j^{\mp}, w_j^{\pm} \right) = 0,
\end{eqnarray}
where $\rho_j = {\rm Re} \left( L_- w_j(x), w_j(x)\right)$.
\end{Lemma}

\begin{Proof}
Since operator $L_-$ is real-valued, the eigenvector $w_j(x)$ of the
problem (\ref{spectrum}) with ${\rm Im}(\lambda_j) > 0$ has both
real and imaginary parts, which are linearly independent. Under the
constraint $\lambda_j \neq 0$,
$$
\left( L_- w_j, \bar{w}_j \right) = \lambda_j^{-1} \left( L_- w_j,
\partial_x L_- \bar{w}_j \right) = 0,
$$
and the orthogonality equations (\ref{orthogonality5}) follow by
direct computations.
\end{Proof}

\begin{Theorem}
\label{proposition-invariance} Let $N_{\rm real}$ be the number of
real eigenvalues of the stability problem (\ref{spectrum}) with
${\rm Re}(\lambda) > 0$, $N_{\rm comp}$ be the number of complex
eigenvalues with ${\rm Re}(\lambda) > 0$ and ${\rm Im}(\lambda)
> 0$, and $N_{\rm imag}^-$ be the number of imaginary eigenvalues
with ${\rm Im}(\lambda) > 0$ and $\rho_j \leq 0$ in
(\ref{orthogonality5}). Assume that all embedded (zero and
imaginary) eigenvalues of the stability problem (\ref{spectrum})
are simple, such that ${\rm Ker}(L_+) = \{\phi\} \in {\cal H}$.
Then,
\begin{equation}
\label{closure-relation1} N_{\rm real} + 2 N_{\rm comp} + 2 N_{\rm
imag}^- = n(L_-) - N_{\rm zero}^-,
\end{equation}
where $N_{\rm zero}^- = 1$ if $(L_-^{-1} \phi,\phi) \leq 0$ and
$N_{\rm zero}^- = 0$ otherwise.
\end{Theorem}

\begin{Proof}
Each isolated and embedded eigenvalue $\gamma_j = - \lambda_j^2$
of the generalized eigenvalue problem
(\ref{generalized_eigenvalue_problem}) is at least double with two
linearly independent eigenvectors $u_j^{\pm}(x)$, such that
$w_j^{\pm} = \partial_x u_j^{\pm}$. By Lemma \ref{lemma-inversion}
and Corollary \ref{remark-complex-KdV}, the dimension of
non-positive invariant subspace of $\Pi_{\kappa}$ for isolated
(real and complex) eigenvalues coincide with the algebraic
multiplicities of isolated eigenvalues, such that $N_n^- = N_p^- =
N_{\rm real}$ and $N_{c^+} = 2 N_{\rm comp}$. By convention in
Remark \ref{embedded-eigenvalues}, all embedded (zero and purely
imaginary) eigenvalues are assumed to be simple. By Lemma
\ref{lemma-orthogonality3} and the relation for eigenvectors of
the stability problem (\ref{spectrum}),
\begin{equation}
\label{relation-L-M-forms} (L_+ u, u) = (L_- u', u') = (L_- w, w),
\end{equation}
we have $N_n^+ = 2 N_{\rm imag}^-$ and $N_n^0 = N_{\rm zero}^0$.
The count (\ref{closure-relation1}) follows by equality
(\ref{negative-index-K}) of Theorem \ref{equality 1}.
\end{Proof}

\begin{Remark}
{\rm Theorem \ref{proposition-invariance} can be generalized to any
KdV-type evolution equation, when the linearization operator $L_-$
is invariant with respect to the transformation $x \mapsto -x$. When
$N_{\rm imag}^- = N_n^0 = 0$, the relation (\ref{closure-relation1})
extends the Morse index theory from gradient dynamical systems to
the KdV-type Hamiltonian systems. For gradient dynamical systems,
all negative eigenvalues of $L_-$ are related to real unstable
eigenvalues of the stability problem. For the KdV-type Hamiltonian
system, negative eigenvalues of $L_-$ may generate both real and
complex unstable eigenvalues in the problem (\ref{spectrum}).}
\end{Remark}

{\bf Example 4.} Let $\phi(x)$ be a single-pulse solution on $x \in
\mathbb{R}$ (which does not have to be positive). Assume that the
operator $L_-$ has a one-dimensional kernel in ${\cal X}$ with the
eigenvector $\phi'(x)$. Then the kernel of $L_+$ is one-dimensional
in ${\cal X}$ with the eigenvector $\phi(x)$, such that $z(L_+) =
1$.

\begin{itemize}
\item Since $(\phi,\phi') = 0$, then ${\rm Ker}(L_+) \in {\cal H}$
and $z_0 = 0$. Since $L_-
\partial_c \phi(x) = - \phi(x)$, we have $z_1 = 1$ if $\frac{d}{dc} \|
\phi \|^2_{L^2} = 0$. Since the latter case violates the assumption
that the embedded kernel is simple in ${\cal H}$, we shall consider
the case $\frac{d}{dc} \| \phi \|^2_{L^2} \neq 0$, such that $z_1 =
0$.

\item Since $(K\phi,\phi) = -(\partial_c \phi,\phi) = -\frac{1}{2}
\frac{d}{dc} \| \phi \|^2_{L^2}$, we have $N_{\rm zero}^- = 1$ if
$\frac{d}{dc} \| \phi \|^2_{L^2} > 0$ and $N_{\rm zero}^- = 0$ if
$\frac{d}{dc} \| \phi \|^2_{L^2} < 0$. By Lemma
\ref{lemma-shift-zero} and Remark \ref{remark-shift-zero}, ${\rm
dim}({\cal H}_{A + \delta K}^-) = {\rm dim}({\cal H}_A^-) + N_{\rm
zero}^-$. By the relation (\ref{negative-index-A}) of Theorem
\ref{equality 1} with $N_p^- = N_{\rm real}$, we have
\begin{equation}
\label{closure-11} N_{\rm real} + 2 N_{\rm comp} + 2 N_{\rm
imag}^- = {\rm dim}({\cal H}_{A + \delta K}^-) - N_{\rm zero}^- =
{\rm dim}({\cal H}_A^-).
\end{equation}

\item  Recall that ${\rm dim}({\cal H}_A^-) = n(L_+) - n_0$ in 
Proposition \ref{lemma-constrained-space} (extended to the case 
$\omega_+ = 0$). By the Sylvester Inertia Law Theorem and the 
relation (\ref{relation-L-M-forms}), we have $n(L_+) = n(L_-)$ in 
${\cal X}$. Since $A(0) = -(L_+^{-1}\phi',\phi') = 
-(L_-^{-1}\phi,\phi)$, we have $n_0 = N_{\rm zero}^-$ in ${\cal H}$. 
As a result, the relation (\ref{closure-11}) recovers the same 
equality (\ref{closure-relation1}).
\end{itemize}

\begin{Remark}
{\rm General stability-instability results for the traveling waves
of the KdV-type equations were obtained in \cite{BSS87,SS90} when
$n(L_-) = 1$. In this case, $N_{\rm comp} = N_{\rm imag}^- = 0$
and $N_{\rm real} = 1 - N_{\rm zero}^-$, such that stability
follows from $\frac{d}{dc} \| \phi \|^2_{L^2} > 0$ and instability
follows from $\frac{d}{dc} \| \phi \|^2_{L^2} < 0$. By a different
method, Lyapunov stability of positive traveling waves $\phi(x)$
was considered in \cite{Wen}. Specific studies of stability for
the fifth-order KdV equation (\ref{general-KdV}) were reported in
\cite{IS92,DK99} with the energy-momentum methods. Extension of
the stability-instability theorems of \cite{BSS87,Wen} with no
assumption on a simple negative eigenvalue of $L_-$ was developed
in \cite{L99,Pava03} with a variational method. The variational
theory is limited however to the case of homogeneous
nonlinearities, e.g. $b_3 = 0$ or $b_1 = b_2 = 0$. Our treatment
of stability in the fifth-order KdV equation (\ref{general-KdV})
is completely new and it exploits a similarity between stability
problems of KdV and NLS solitons. The pioneer application of the
new theory to stability of N-solitons in the KdV hierarchy is
reported in \cite{KP}. Further progress on the same topic will
appear soon in \cite{Dias,Sand}. }
\end{Remark}

{\bf Acknowledgement:} This work was completed with the support of
the SharcNet Graduate Scholarship and the PREA grant. The authors
are grateful to Prof. T. Azizov for his help in the proof of the
Pontryagin's Theorem.

\end{document}